\RequirePackage{etoolbox}
\csdef{input@path}{{style/}{graphics/}}
\documentclass[aap,MSNbibl,dvips]{arximspdf}
\makeatletter
   \@ifpackageloaded{graphicx}{}{\usepackage{graphicx}}
\makeatother
\usepackage{accents}%natbib

\doi{10.1214/14-AAP1032} %kopijuoti is PTS
\volume{25}
\issue{3}
\pubyear{2015}
\firstpage{1616}
\lastpage{1649}

\makeatletter
\newtheorem{theorem}{Theorem}%[section]
\newtheorem{proposition}[theorem]{Proposition}
\newtheorem{lemma}[theorem]{Lemma}
\newtheorem{corollary}[theorem]{Corollary}
\newproclaim{rem}{Remark}
\newcommand{\rvec}[1]{\accentset{\rightarrow}{#1}}
\newcommand{\lvec}[1]{\accentset{\leftarrow}{#1}}
\newcommand{\var}{\operatorname{Var}}
\makeatother

\begin{document}
\begin{frontmatter}

\title{A particle system with cooperative branching and~coalescence\thanksref{T1}}
\runtitle{Cooperative branching-coalescent}

\begin{aug}
\author[A]{\fnms{Anja}~\snm{Sturm}\corref{}\ead[label=e1]{asturm@math.uni-goettingen.de}}
\and
\author[B]{\fnms{Jan M.}~\snm{Swart}\ead[label=e2]{swart@utia.cas.cz}\thanksref{T2}}
\runauthor{A. Sturm and J. M. Swart}
\affiliation{Georg-August-Universit\"at G\"ottingen and Institute of
Information Theory and~Automation of the ASCR (\'UTIA)}
\address[A]{Institute for Mathematical Stochastics\\
Georg-August-Universit\"at G\"ottingen\\
Goldschmidtstr.~7\\
37077 G\"ottingen\\
Germany\\
\printead{e1}} %adresu isvedimo komanda gale!
\address[B]{Institute of Information Theory\\
\quad and Automation of the ASCR (\'UTIA)\\
Pod vod\'arenskou v\v e\v z\' i 4\\
18208 Praha 8\\
Czech Republic\\
\printead{e2}}
\end{aug}
\thankstext{T1}{Supported in part by DFG priority programme 1590.}
\thankstext{T2}{Work sponsored by GA\v CR Grant P201/10/0752.}

% HISTORY:
\received{\smonth{11} \syear{2013}}
\revised{\smonth{4} \syear{2014}}

% ABSTRACT
%
\begin{abstract}
In this paper, we introduce a one-dimensional model of particles performing
independent random walks, where only pairs of particles can produce offspring
(``cooperative branching''), and particles that land on an \mbox{occupied}
site merge
with the particle present on that site (``coalescence''). We show that the
system undergoes a phase transition as the branching rate is increased.
For small branching rates, the upper invariant law is trivial, and the process started
with finitely many particles a.s. ends up with a single particle.
Both statements are not true for high branching rates. An interesting
feature of the process is that the spectral gap is zero even for low branching
rates. Indeed, if the branching rate is small enough, then we show
that for the process started in the fully occupied state, the particle density
decays as one over the square root of time, and the same is true for
the decay of the
probability that the process still has more than
one particle at a later time if it started with two particles.
\end{abstract}

% KEYWORDS
\end{frontmatter}

\setcounter{footnote}{2}

%s1 #&#
\section{Introduction and main results}

%s1.1 #&#
\subsection{Definition of the model}\label{Sdefs}

\begin{figure}[b]

\includegraphics{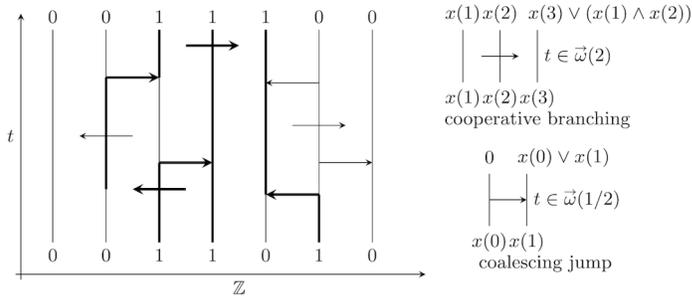}

%
%
%
%
%
%
%
%
%
%
%
%x(2))$};
%
%
%
%
\caption{Example of a graphical representation with explanation of the rules.
Bold lines indicate the presence of a particle. Arrows that are used by a
particle to jump through, or by a pair of particles to give birth to a third
particle are also drawn bold, regardless of whether such a particle
lands on
an occupied site or not.}
\label{figgraph}
\end{figure}

Let $\{0,1\}^{\mathbb Z}$ be the space of all configurations $\ldots
10010101101\ldots$
of zeros and ones on the integers. We denote such a configuration by
$x=(x(i))_{i\in{\mathbb Z}}$ with $x(i)\in\{0,1\}$. Let $\lambda\geq
0$ be a parameter, to
be referred to as the \emph{cooperative branching rate}. We will be
interested in the continuous-time Markov process $X=(X_t)_{t\geq0}$ taking
values in $\{0,1\}^{\mathbb Z}$ and with right-continuous sample paths,
such that if
$X$ is in the state $x$, then for each $i\in{\mathbb Z}$, it makes
transitions with the following exponential rates:
%
%e1 #&#
\begin{equation}\label{ratedef} %
\begin{array} {rll}
\mbox{if }x(i)=1,\qquad\mbox{then}&
\bigl(x(i),x(i+1) \bigr)\mapsto(0,1) \qquad&\mbox{at rate }\frac{1}{2},
\\[7pt]
& \bigl(x(i-1),x(i) \bigr)\mapsto(1,0) \qquad&\mbox{at rate }\frac{1}{2},
\\[7pt]
\mbox{if } \bigl(x(i),x(i+1) \bigr)=(1,1),\qquad\mbox{then}&
x(i+2)\mapsto1
\quad&\mbox{at rate }\frac{1}{2}\lambda,
\\[7pt]
&x(i-1)\mapsto1 \quad&\mbox{at rate }\frac{1}{2}\lambda.
\end{array}\hspace*{-18pt}
\end{equation}
In these transitions, $x(j)$ remains the same for all sites $j$ not listed.
We may construct such a process with the help of a graphical
representation as
follows. For each $i\in{\mathbb Z}/2=\{k/2\dvtx k\in{\mathbb Z}\}$, let
$\rvec\omega(i),\lvec\omega(i)\subset{\mathbb R}$ be Poisson
subsets of the real line. We assume
that all these Poisson sets are independent and that $\rvec\omega
(i),\lvec\omega(i)$
have\vspace*{1pt} intensity $\frac{1}{2}$ if $i\in{\mathbb Z}+\frac{1}{2}:=\{
k+\frac{1}{2}\dvtx k\in{\mathbb Z}\}$ and
intensity $\frac{1}{2}\lambda$ if $i\in{\mathbb Z}$. In pictures, we
plot ${\mathbb Z}$ horizontally and
time vertically. We indicate the presence of a point $t\in\rvec\omega(i)$
[resp., $t\in\lvec\omega(i)$] by drawing a vector at time $t$ from
$i-\frac{1}{2}$ to
$i+\frac{1}{2}$ (resp., from $i+\frac{1}{2}$ to $i-\frac{1}{2}$);
see Figure~\ref{figgraph}.

%
%f1 #&#

We interpret the points of $\rvec\omega(i),\lvec\omega(i)$ with
$i\in{\mathbb Z}+\frac{1}{2}$ as \emph{coalescing jump events} and
those with $i\in{\mathbb Z}$
as \emph{cooperative branching events}. Starting from an initial state
$X_0=x\in\{0,1\}^{\mathbb Z}$ at time zero, we construct a process
$X^x=X=(X_t)_{t\geq
0}$ that changes its state only at coalescing jump events and cooperative
branching events according to the following rules.

If immediately prior to some coalescing jump event $t\in\rvec\omega
(i)$ (with
$i\in{\mathbb Z}+\frac{1}{2}$) the state is $X_{t-}$ and
$X_{t-}(i-\frac{1}{2})=1$,
then we set $X_t(i-\frac{1}{2})=0$, $X_t(i+\frac{1}{2})=1$.
Everywhere else,
we do nothing; that is, we set $X_t(j)=X_{t-}(j)$ for all \mbox{$j\neq
i-\frac{1}{2}$}, $i+\frac{1}{2}$. If
$X_{t-}(i-\frac{1}{2})=0$, then we set $X_t(j)=X_{t-}(j)$ for all $j$;
that is,
we do nothing. Interpreting a one (resp., zero) as the presence
(resp., absence) of a particle, this says that at each time $t\in\rvec
\omega(i)$,
any particle that may be present at the site $i-\frac{1}{2}$ jumps to
$i+\frac{1}{2}$, coalescing with any particle that may already be present
there. Likewise, at times $t\in\lvec\omega(i)$ a particle at
$i+\frac{1}{2}$ (if
there is one)\vspace*{1pt} jumps to $i-\frac{1}{2}$.

If immediately prior to some cooperative branching event $t\in\rvec
\omega(i)$
(with $i\in{\mathbb Z}$) we have $(X_{t-}(i-1),X_{t-}(i))=(1,1)$, then
we set
$X_t(i+1)=1$ and $X_t(j)=X_{t-}(j)$ for all $j \neq i+1$. If
$(X_{t-}(i-1),X_{t-}(i))\neq(1,1)$, then we do nothing. We may also describe
this by saying that if $i-1$ and $i$ are both occupied by a particle, then
these two particles cooperate to produce a particle at $i+1$, which coalesces
with any particle that may already be present there. Likewise, at times
$t\in\lvec\omega(i)$, if there are particles at both $i$ and $i+1$,
then these
give birth to a particle at $i-1$.

These rules are further illustrated in Figure~\ref{figgraph},
together with
an example of a graphical representation. It can be checked by standard means\footnote{Essentially, one can check that for given $s\leq t$, the
number of
sites $j$ whose state at time $s$ could possibly influence the state of a
given site $i$ at time $t$ is a.s. finite, and in fact its expectation grows
at most exponentially in $t-s$.} that the graphical representation yields,
for each initial state $x\in\{0,1\}^{\mathbb Z}$, a well-defined \mbox{$\{
0,1\}^{\mathbb Z}$-}valued
Markov process $X^x=(X^x_t)_{t\geq0}$ with initial state $X^x_0=x$.
Note that
the graphical representation provides a natural coupling between processes
started in different (deterministic) initial states. The graphical
representation can also be used to construct processes started in random
initial states. In this case the initial state must be independent of the
graphical representation. We call our process the \emph{cooperative
branching-coalescent} with \emph{cooperative branching rate $\lambda
$}. Our
motivation for studying this particular model will be explained in
detail in
Section~\ref{Smotiv} below.

It\vspace*{2pt} will often be convenient to use set notation for our state
space. Identifying a set $A\subset{\mathbb Z}$ with its indicator
function $1_A$, we may
identify the space $\{0,1\}^{\mathbb Z}$ with the space ${\mathcal
P}({\mathbb Z})$ of all subsets of
${\mathbb Z}$. For each $A\subset{\mathbb Z}$, we let
%
%e2 #&#
\begin{equation}
\label{etadef} \eta^A_t:=\bigl\{i\in{\mathbb
Z}\dvtx X^{1_A}_t(i)=1\bigr\}\qquad(t\geq0)
\end{equation}
denote the set of occupied sites at time $t$ for the process started
with the
initial set of occupied sites $A$. Then $\eta^A=(\eta^A_t)_{t\geq0}$
is just a\vspace*{1pt}
different notation for the cooperative branching-coalescent
$X^{1_A}=(X^{1_A}_t)_{t\geq0}$. Because of certain notational
advantages, we
will usually (but not always) use this sort of set notation for our processes.

%Here $\wedge$ and $\vee$ denote the minimum and maximum, respectively.
%Coalescing jump arrows or cooperative branching arrows that are used
%(in the
%sense that a particle jumps or a pair of particles gives birth to a
%third
%particle

%s1.2 #&#
\subsection{Basic facts}

Recall \cite{Lig85}, Theorem~II.2.4, that the laws $\mu:={\mathbb
P}[Y\in\cdot]$ and
$\nu:={\mathbb P}[Z\in\cdot]$ of two $\{0,1\}^{\mathbb Z}$-valued
random variables $Y$ and $Z$
are said to be stochastically ordered, denoted as $\mu\leq\nu$, if
$Y$ and $Z$
can be coupled such that $Y \leq Z$ a.s., by which we mean that
$Y(i)\leq Z(i)$ $(i\in{\mathbb Z})$ a.s. Equivalently, using
set notation, this says that the laws of two ${\mathcal P}({\mathbb
Z})$-valued random
variables $\eta,\xi$ are stochastically ordered if they can be
coupled such
that $\eta\subset\xi$. It is a simple consequence of our graphical
representation
that cooperative branching-coalescents are monotone in the following sense.

%le1 #&#
\begin{lemma}[(Monotonicity)]\label{Lmonot}
Let $\eta$ and $\eta'$ be cooperative branching-coales\-cents
with cooperative
branching rates $\lambda$ and $\lambda'$, respectively. Assume that
$\lambda\leq\lambda'$
and ${\mathbb P}[\eta_0\in\cdot]\leq{\mathbb P}[\eta'_0\in\cdot
]$. Then
${\mathbb P}[\eta_t\in\cdot]\leq{\mathbb P}[\eta'_t\in\cdot
]$ for all $t\geq0$.
\end{lemma}

\begin{pf}
We first use the fact that
${\mathbb P}[\eta_0\in\cdot]\leq{\mathbb P}[\eta'_0\in\cdot
]$ to couple $\eta_0$ and
$\eta'_0$ in such a way that $\eta_0\subset\eta'_0$ a.s. Next,\vspace*{2pt} we
construct a
graphical representation, consisting of Poisson sets $\rvec\omega
(i),\lvec\omega(i)$
and $\rvec\omega'(i),\lvec\omega'(i)$, respectively, for the
processes $\eta$ and
$\eta'$, independent of $(\eta_0,\eta'_0)$, in the following way. Starting
from a graphical representation for $\eta$, we define
$\rvec\omega'(i):=\rvec\omega(i)$ and $\lvec\omega'(i):=\lvec
\omega(i)$ for $i\in{\mathbb Z}+\frac{1}{2}$,
that is, the processes $\eta$ and $\eta'$ use the same coalescing jump
events. For $i\in{\mathbb Z}$, we let $\rvec\omega''(i)$ and $\lvec
\omega''(i)$ be
independent Poisson sets with intensity $\frac{1}{2}(\lambda'-\lambda
)$ and set
$\rvec\omega'(i):=\rvec\omega(i)+\rvec\omega''(i)$ and likewise
$\lvec\omega'(i):=\lvec\omega(i)+\lvec\omega''(i)$. In this way,
the cooperative
branching events of $\eta$ are a subset of those of $\eta'$. It is now
straightforward to check from the rules of a graphical representation that
$\eta_t\subset\eta'_t$ a.s. for each $t\geq0$.
\end{pf}

It is easy to check that the rules of our graphical representation moreover
imply the following property.

%le2 #&#
\begin{lemma}[(Subadditivity)]\label{Lsubad}
For a given graphical representation, one has
%
%e3 #&#
\begin{equation}
\label{subad} \eta^A_t\cup\eta^B_t
\subset\eta^{A\cup B}_t\qquad(t\geq0, A,B\subset{\mathbb Z}).
\end{equation}
\end{lemma}

Processes that have a graphical representation for which equality holds in
(\ref{subad}) are called \emph{additive} \cite{Gri79},
Proposition~II.1.2. Our
process, however, only has the weaker property (\ref{subad}) (unless
$\lambda=0$
which is a pure coalescing random walk). It can, moreover, be checked that
because of the coalescing random walk dynamics, which involves jumps between
incomparable states, meaning that jumps occur from state $x$ to $x'$
such that
neither $x \leq x'$ nor $x' \leq x$, our process does not satisfy~\cite{Lig85}, formula~(II.2.19), and hence it does not preserve positive
correlations.

Lemma~\ref{Lmonot} with $\lambda=\lambda'$ says that the cooperative
branching-coalescent is a monotone interacting particle system. It is
well known that this implies the existence of an invariant law
$\bar\nu$, called the \emph{upper invariant law}, such that
%
%e4 #&#
\begin{equation}
\label{upper} {\mathbb P} \bigl[\eta^{\mathbb Z}_t\in\cdot\bigr]
\mathop{\Longrightarrow}_{{t}\to\infty}\bar\nu.
\end{equation}
(For attractive spin systems, this is proved in
\cite{Lig85}, Theorem~II.2.3. Although not stated there, the proof actually
carries over without a change to any monotone interacting particle system.)
Here, $\Rightarrow$ denotes weak convergence of probability measures on
$\{0,1\}^{\mathbb Z}$, equipped with the product topology. Moreover,
$\bar\nu$ dominates
any other invariant law of the process in the stochastic order (hence its
name). Using again Lemma~\ref{Lmonot}, but now with $\lambda\leq
\lambda'$, it is,
moreover, easy to see that the upper invariant laws
$\bar\nu_\lambda,\bar\nu_{\lambda'}$ corresponding to
cooperative branching rates
$\lambda\leq\lambda'$ are stochastically ordered as
$\bar\nu_\lambda\leq\bar\nu_{\lambda'}$. We say that
$\bar\nu$ is \emph{nontrivial} if
$\bar\nu$ gives zero probability to the empty configuration,
that is, if
$\bar\nu(\{\varnothing\})=0$, and we let
%
%e5 #&#
\begin{equation}
\label{dens} \theta(\lambda):=\int\bar\nu_\lambda(\mathrm{d}A)1_{\{0\in
A\}}
\end{equation}
denote the probability under $\bar\nu$ of finding a particle in
the origin.

It is clear from our dynamics that a process started with a single particle
will consist of a single particle at all times, and this particle performs
simple random walk on ${\mathbb Z}$. We will say that the process \emph
{survives}
for a given value $\lambda$ of the cooperative branching rate if the
probability
%
%e6 #&#
\begin{equation}
\label{surviv} \psi(\lambda):={\mathbb P} \bigl[\bigl|\eta^{\{0,1\}}_t\bigr|
\geq2\ \forall t\geq0 \bigr]
\end{equation}
is positive. If the process does not survive, then we say that it \emph{dies
out}. (Even though, of course, there will always be one particle left. But
since only pairs of particles can branch or coalesce, we are naturally
interested in whether there will always survive at least two particles
in the system.) It is easy to see from Lemma~\ref{Lmonot} that this
probability is nondecreasing in the cooperative branching rate $\lambda$.

%s1.3 #&#
\subsection{Main results}

Our first main result says that the cooperative bran\-ching-coalescent exhibits
a phase transition, both in terms of its upper invariant law and in
terms of
survival.

%th3 #&#
\begin{theorem}[(Phase transition)]\label{Tphase}
\textup{(a)} There exists a $1\leq\lambda_{\mathrm{c}}<\infty$
such that
$\bar\nu_\lambda=\delta_\varnothing$ for $\lambda<\lambda
_{\mathrm{c}}$, but $\bar\nu_\lambda$ is
nontrivial for $\lambda>\lambda_{\mathrm{c}}$.

\textup{(b)} There exists a $1\leq\lambda'_{\mathrm{c}}<\infty$
such that the
process dies out for $\lambda<\lambda'_{\mathrm{c}}$ and survives for
$\lambda>\lambda'_{\mathrm{c}}$.
\end{theorem}

The basic idea behind the proof of Theorem~\ref{Tphase}, which can be
found in
Section~\ref{Sphasetrans}, is easily
explained. If $\lambda<1$, then each pair of particles on neighboring positions
on average creates fewer particles by cooperative branching than are
lost by
coalescence, from which it is not too hard to conclude that no nontrivial
invariant law is possible, and systems started with finitely many
particles end
up with one particle a.s.; see Section~\ref{Sextinction}.
On the other hand, for sufficiently high cooperative branching
rates, a pair of particles on neighboring positions has a high
probability of
producing particles on both of its neighboring sites before any of its
particles makes a jump. Using this, one can set up a comparison with
supercritical oriented percolation which gives both survival and
existence of
a nontrivial invariant law. This is done in Section~\ref{Ssurvival}
where we also complete the proof
of Theorem~\ref{Tphase}.

\begin{figure}%[t]

\includegraphics{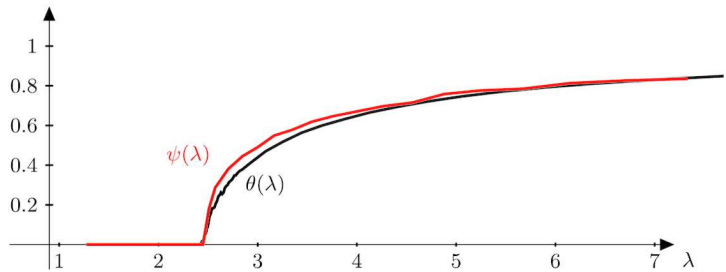}

%%\draw(1,1.2) node[right] {$\tet(\la)$};
%%\draw(1.649,0) node[below] {$\la_{\mathrm{c}}$};
%%\draw[very thick] (-0.07,1) -- (0.07,1);
%%\draw(0,1) node[left] {$1$};
%%\draw(1,0) node[below right] {1};
\caption{Density $\theta(\lambda)$ of the upper invariant law and survival
probability $\psi(\lambda)$ (plotted in black and red, resp.) of the
cooperative branching-coalescent as a function of the cooperative branching
rate.}
\label{fignumer}
\end{figure}

We do not know if $\lambda_{\mathrm{c}}=\lambda'_{\mathrm{c}}$, although it
seems plausible
that this is indeed the case. Numerically, both critical points are
given by
%
%e7 #&#
\begin{equation}
\lambda_{\mathrm{c}}\approx\lambda'_{\mathrm{c}}\approx2.47
\pm0.02;
\end{equation}
see Figure~\ref{fignumer}.
%
%f2 #&#

\begin{figure}[t]

\includegraphics{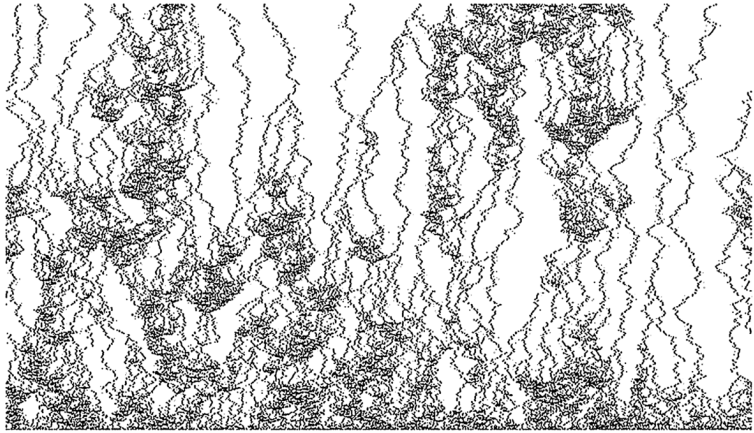}

\caption{Simulation of a near-critical cooperative branching-coalescent with
$\lambda=7/3$ on a lattice of 700 sites with periodic boundary
conditions, started from the fully occupied initial state. Space is plotted
horizontally, time vertically and black indicates the presence of a
particle.}
\label{figdraw}
\end{figure}

Superficially, the behavior of the cooperative branching-coalescent looks
similar to that of the contact process, but the critical exponent associated
with the density of the upper invariant law seems to be different. For the
one-dimensional contact process, and indeed for many other, similar particle
systems that are supposed to be in the same universality class, it is believed
(and explained by nonrigorous renormalization group theory) that the density
of the upper invariant law grows like $(\lambda-\lambda_{\mathrm{c}})^\beta$ with
$\beta\approx0.27648$ \cite{Hin00}, Section~3.4. For the cooperative
branching-coalescent, this critical exponent $\beta$ [as read off from
a plot of
$\log\theta$ versus $\log(\lambda-\lambda_{\mathrm{c}})$] seems to be
approximately
$\beta\approx0.5\pm0.1$. A picture of a near-critical process is
shown in
Figure~\ref{figdraw}.

%
%f3 #&#

From a physical point of view, different critical exponents are to be
expected since we can prove that our process exhibits critical behavior (in
particular, the absence of a spectral gap) in the whole regime
$\lambda\leq\lambda_{\mathrm{c}}$. This is our second main result which
implies that
the behavior of the cooperative branching-coalescent
contrasts with the behavior of the contact process that
is known to have a spectral gap in the whole subcritical regime.
Indeed, the probability that a subcritical contact process
started with finitely many infected sites survives until time $t$ decays
exponentially in $t$ \cite{BG91}, and by the self-duality of the contact
process, the same is true for the density at time $t$ of the process started
with all sites occupied. Our result shows that for the cooperative
branching-coalescent, both quantities decay according to a power law with
exponent $-1/2$.

%th4 #&#
\begin{theorem}[(Decay rate in the subcritical regime)]\label{Tdecay}
Let $\eta^{\{0,1\}}$ and $\eta^{\mathbb Z}$ be cooperative branching-coalescents with
cooperative branching rate $\lambda\geq0$, started with two particles at
neighboring sites or in the fully occupied state, respectively. Then there
exists a constant $c>0$ such that for all $\lambda\geq0$,
%
%e8 #&#
\begin{equation}
\label{lbd} {\mathbb P} \bigl[\bigl|\eta^{\{0,1\}}_t\bigr|\geq2 \bigr]
\geq ct^{-1/2} \quad\mbox{and}\quad{\mathbb P}\bigl[0\in
\eta^{\mathbb Z}_t\bigr]\geq ct^{-1/2}\qquad(t\geq0).
\end{equation}
Moreover, there exists a constant $C<\infty$ such that for each $0\leq
\lambda<
1/2$,
%
%e9 #&#
\begin{equation}
\label{ubd} {\mathbb P} \bigl[\bigl|\eta^{\{0,1\}}_t\bigr|\geq2 \bigr]
\leq Ct^{-1/2} \quad\mbox{and}\quad{\mathbb P}\bigl[0\in
\eta^{\mathbb Z}_t\bigr]\leq Ct^{-1/2}\qquad(t\geq0).
\end{equation}
\end{theorem}

The proof of (\ref{lbd}) is easy: by Lemma~\ref{Lmonot}, we can estimate
$\eta^{\mathbb Z}$ from below by a system with cooperative branching
rate zero, that is,
by a pure coalescencing random walk, for which the decay of both
quantities is well known
to follow a power law with exponent $-1/2$.
The proof of (\ref{ubd}) is more involved and depends on estimating the
survival probability of a somewhat complicated ``superdual'' process.
The proof of Theorem~\ref{Tdecay}
is completed at the end of Section~\ref{Sdecay}.

%s1.4 #&#
\subsection{Discussion and motivation}\label{Smotiv}

Systems with cooperative branching, but different death mechanisms have been
considered before in the literature. In particular, Noble \cite{Nob92} and
Neuhauser \cite{Neu94} have studied a ``sexual reproduction process''
in which
particles perform cooperative branching (but no coalescence) and die with
constant rate. The name of this process is a bit misleading since it
does not
distinguish organisms with different sexes. An interesting feature of
it is
that the corresponding mean-field model exhibits a first order phase
transition, which is reflected in meta-stable behavior of the spatial model
with strong mixing or long-range interaction.

In the physics literature, considerable attention has been paid to the ``pair
contact process'' where again the reproduction mechanisms is cooperative
branching, but the death mechanism is annihilation (two particles are
simultaneously removed). Whether this model belongs to the directed
percolation (DP) universality class is the subject of ongoing debate
\cite{SB08,Par12}.

Our motivation for studying the cooperative branching-coalescent is
multifaceted. As detailed below, we regard the model as an
interesting toy model in and by itself, both from a biological
and mathematical perspective. In addition, the model
is of relevance due to connections to other interesting models for
which it is
potentially harder to prove the results that we can obtain here.

From a purely mathematical perspective, the cooperative branching-\break coalescent
is interesting because of the critical behavior in the extinction phase
(proved in Theorem~\ref{Tdecay}), which sets it apart from more usual models
exhibiting a phase transition between extinction and survival, such as the
contact process. This criticality arises from the fact that extinction is
driven by coalescence alone, a property presumably shared with other models
that are more difficult to treat.

From a biological perspective, we have two, rather different
motivations for
studying cooperative branching. First of all, interpreting particles as
organisms, we may view the cooperative branching-coalescent as a model for
population dynamics. The assumption that only pairs of individuals can
reproduce is, of course, rather natural. Although usually, the members
of such
a pair need to be of opposite sex (a fact not incorporated in our model),
there are in fact quite a lot of organisms (such as snails) that are
hermaphroditic, that is, each individual plays the role of both sexes,
but that do
not self-fertilize. In this interpretation, the random walk dynamics models
dispersal of organisms while the coalescence represents a death rate
that is
quadratic in the local population size. Such a quadratic death rate naturally
models deaths due to competition between individuals for limited space and
resources \cite{AS05}.

To make the model more realistic, one could (as in \cite{AS05}) also
add a
linear component in the death rate, representing spontaneous deaths
that are
not due to competition with other individuals. Doing so would, however,
radically change the properties of the model. In particular, this would
destroy the validity of Theorem~\ref{Tdecay} and presumably yield a
model in
the universality class of the contact process. For these reasons, we have
chosen not to do this.

Our second biological motivation for considering the cooperative
bran\-ching-coalescent comes from the study of balancing selection (sometimes
also called heterozygosity selection or negative frequency dependent
selection). This is the \mbox{phenomenon} that genetically similar individuals often
compete more strongly with each other than with genetically more different
individuals. This could, for example, be due to the fact that
genetically more
different individuals need a somewhat different set of resources for survival.

In order to model this effect, Neuhauser and Pacala \cite{NP99}
introduced a
variation of the voter model in which types that are locally in the minority
have an advantage (due to the presumed smaller competition with
neighbors). A
very similar model, dubbed the ``rebellious voter model,'' was
introduced in
\cite{SS08}. Numerical simulations backed up, in part, by rigorous
mathematics (see \cite{SV10} and references therein) have shown that
typically, such models in dimensions $d\geq2$ tend to have an
invariant law
in which both types are present for all values of the selection
parameter, but
in dimension one undergo a phase transition between noncoexistence and
coexistence as the selective advantage for locally rare types is increased.

Proving the existence of this phase transition, and in particular the
existence of the noncoexisting phase, has proved to be difficult,
however. Both through duality and by considering the corresponding ``interface
model'' (as explained in \cite{SS08}), noncoexistence can be shown to be
equivalent to the extinction of a branching--annihilating particle system,
where single particles give birth to two offspring at once, and pairs of
particles annihilate each other, with certain rates. Such systems are
parity-preserving (i.e., even/oddness of the initial number of
particles is
conserved), so extinction needs to be interpreted as starting from even
initial states, since odd systems can never die out completely. These systems
are similar to the cooperative branching-coalescent in the sense that single
particles cannot die, and hence extinction relies on the recurrence of
one-dimensional random walk. It seems that in the extinction regime, these
systems effectively behave like a small perturbation of annihilating random
walks without branching. In particular, it seems likely that their
density and
survival probability (started from an even number of particles) decay like
$t^{-1/2}$, just as for the cooperative branching-coalescent
(Theorem~\ref{Tdecay}).

Contrary to the cooperative branching-coalescent, however, none of these
statements are proven. This is mainly due to two difficulties. First, these
parity-preserving branching--annihilating particle systems are not
monotone, so
the usual coupling arguments fail, and in general one does not even
know if
increasing the branching rate makes survival more likely (although this is
certainly what one sees in all simulations). Second, in a parity-preserving
branching--annihilating particle system, single particles can still branch,
even though most of the particles created in such branchings are
believed to
be quickly lost again due to annihilation. This is related to the
problem of
(strong) interface tightness for the rebellious voter model, which has
recently been shown to imply noncoexistence \cite{Swa13} (although it remains
an open problem to show either occurs).

Using the results of our present paper, we can describe a simple
variation of a
one-dimensional voter model in which rare types have an advantage and for
which the existence of a phase transition between noncoexistence and
coexistence can be proved.

Consider a one-dimensional, nearest neighbor multitype voter model in which
initially each site has a different type. We assume the usual voter model
dynamics, that is, the type of each site is updated with rate one, at
which event
it is replaced by the type found on either side immediately to the left or
right of it, with equal probabilities. In addition, with rate $\lambda
$, we
assume that each singleton, that is, each site that is occupied by a
type that
occurs nowhere else, gives birth to a completely new type which is placed
on one of its neighboring sites, with equal probabilities. Let $Y_t(i)$ denote
the type of site $i$ at time $t$ in this model. Then a little thinking
convinces one that
%
%e10 #&#
\begin{equation}
\label{inter} \eta^{\mathbb Z}_t:=\bigl\{i\in{\mathbb
Z}\dvtx Y_t(i)\neq Y_t(i+1)\bigr\}\qquad(t\geq0)
\end{equation}
defines a cooperative branching-coalescent with cooperative branching rate
$\lambda$, started in the fully occupied state. Note
that $\eta^{\mathbb Z}_t$ is the set of interfaces of $Y_t$, that is,
boundaries where
different types meet. Now Theorem~\ref{Tphase}(a) together with
(\ref{upper}) show that $Y_t$ tends in law to a constant configuration for
$\lambda<\lambda_{\mathrm{c}}$ but to an invariant law in which different
types coexist
for $\lambda>\lambda_{\mathrm{c}}$.

This model is obviously somewhat artificial since it depends crucially
on the
nearest-neighbor property of the interaction, which implies that at all times
each type present in the population occupies a single interval.
Moreover, the
assumption that singletons give birth to a \emph{new} type is not well
motivated from the biological point of view. Nevertheless, the general
behavior of the model seems to be similar to that of other, better-motivated
models with balancing selection such as the rebellious voter model of
\cite{SS08}. In fact, we may view the model we have just described as a
variation on the rebellious voter model in which interface tightness and
monotonicity of the interface model have been built in artificially. As such,
we hope that it may also shed some light on this and similar models.

In this context we should also mention another related one-dimensional model,
the cooperative caring double branching annihilating random walk (ccDBARW),
that was recently introduced and analyzed by Blath and Kurt in \cite
{BK11} and
that partially motivated our present paper. In this model, new
particles can
only be created by clumps of at least two particles at neighboring
sites. In
contrast to our cooperative branching-coalescent they are created in pairs
(double branching) on either side of the clump. In addition, particles perform
a random walk. Unlike in our model, particles that land on the same
site do not
coalesce but annihilate each other. The dynamics of the ccDBARW are somewhat
complicated and contrary to the cooperative branching-coalescent it
cannot be
started in infinite initial states. One of its motivations is that it
demonstrates rather dramatically the nonmonotonicity of the classical DBARW,
another parity-preserving branching--annihilating system. Blath and Kurt showed
that the ccDBARW has parameter ranges for survival as well as for extinction,
which implies that at least one phase transition between survival and
extinction must occur. (However, note that due to lack of monotonicity
in this
model a scenario with multiple phase transitions cannot be ruled out.)

%s1.5 #&#
\subsection{Open problems}

The cooperative branching-coalescent has certain nice properties, such as
monotonicity (see Lemma~\ref{Lmonot}), which allow us to give rather short
proofs of Theorems~\ref{Tphase} and \ref{Tdecay}(a). Beyond these basic
facts, however, many questions concerning the model remain open and
seem to
require a substantially bigger effort to be solved. One of the
difficulties of
the cooperative branching-coalescent is the lack of a simple dual
model, such
as one has for the contact process (which is self-dual) or for the rebellious
voter model \cite{SS08}. A result by Gray \cite{Gra86} that holds for general
attractive spin systems implies that the cooperative branching-coalescent has
a dual taking values in the set of finite collections of finite subsets of
${\mathbb Z}$, but this is a fairly complicated process to work with.
In the present
paper, we content ourselves with a process that is only a ``subdual'' (as
explained in Section~\ref{Stermin} below), which nonetheless shows
that some
properties of the full dual can be controlled and which also provides the
basis of the proof of our Theorem~\ref{Tdecay}(b). We intend to
discuss the
relation of our subduality with Gray's full dual in a separate paper. Much
progress in the understanding of cooperative branching-coalescents can be
expected from a better understanding of the full dual process.

In this section, we list and discuss a number of open problems
concerning the
cooperative branching-coalescent.
\begin{longlist}
\item[(P1)] Generalize Theorems~\ref{Tphase} and \ref
{Tdecay} to higher
dimension.
\end{longlist}
It is not hard to define generalizations of the cooperative
branching-coalescent to higher-dimensional lattices ${\mathbb
Z}^d$ $(d\geq2)$. For
such models, the basic Lemmas~\ref{Lmonot} and~\ref{Lsubad} will remain
true, and also Theorem~\ref{Tphase} can probably be generalized
without too much
difficulty. In transient dimensions $d\geq3$, a bit of care is needed in
defining survival, since it is possible that two or more particles separate
forever. The right definition of survival now seems to be that with positive
probability there are pairs of particles at neighboring positions at arbitrary
late times. Generalizing Theorem~\ref{Tdecay} to higher dimensions is less
straightforward since even for the pure coalescent it is known that the decay
of the density has a different asymptotics now, namely $t^{-1} \log t$
in $d=2$ and
$t^{-1}$ in dimensions $d\geq3$; see \cite{BG80}.

\begin{longlist}
\item[(P2)] Prove equality $\lambda_{\mathrm{c}}=\lambda'_{\mathrm{c}}$ of the critical
parameters from Theorem~\ref{Tphase}.
\end{longlist}
Even an inequality in either way would be interesting here. For the contact
process, the analogous result is a simple consequence of self-duality,
which is
not available here. One possible approach is through the following problem.

\begin{longlist}
\item[(P3)] Prove that survival implies a positive edge speed.
\end{longlist}
Here, a positive edge speed means that for the process started with
only the
negative axis occupied,
%
%e11 #&#
\begin{equation}
\liminf_{t\to\infty}t^{-1}\sup\bigl(\eta^{-{\mathbb N}}_t
\bigr)>0.
\end{equation}
This sort of a result could potentially be used to set up a comparison with
supercritical oriented percolation. This is related to the work of
Bezuidenhout, Gray and Grimmett \cite{BG90} and \cite{BG94}, which, however,
does not easily generalize to our model because of the lack of positive
correlations. A more modest problem is whether $\lambda>\lambda_{\mathrm{c}}$ or
$\lambda>\lambda'_{\mathrm{c}}$ (or both) imply a positive edge speed.

\begin{longlist}
\item[(P4)] Prove any estimate for the critical exponent
associated with
the density of the upper invariant law or the survival probability.
\end{longlist}

This looks like a hard problem but any argument that allows one to compare
with the contact process (believed $\beta\approx0.27648$) or
rebellious voter
model (conjectured $\beta\approx0.9$--$1.0$, see the discussion in
\cite{SV10}) would be valuable.

\begin{longlist}
\item[(P5)] For $\lambda>\lambda_{\mathrm{c}}$, show that
$\bar\nu$ is the only
nontrivial translation invariant stationary law, and the limit law started
from any nontrivial translation invariant initial law.
\end{longlist}
This can usually be proved provided one has sufficient control on the dual
model see, for example, the classical proof for the contact process
\cite{Dur80,DG82} or Theorem~5 of
\cite{SS08} for the rebellious voter model. For sufficiently large
$\lambda$, a simpler
proof may be available using monotonicity.

\begin{longlist}
\item[(P6)] Extend the statements in Theorem~\ref
{Tdecay}(b) to all
$\lambda<\lambda'_{\mathrm{c}}$, respectively, $\lambda<\lambda_{\mathrm{c}}$.
\end{longlist}

Again, good control of the dual seems key here.

%s2 #&#
\section{Proof of the phase transition}\label{Sphasetrans}

In this section we prove Theorem~\ref{Tphase} by first showing extinction
(resp., the triviality of the upper invariant law) for small $\lambda$ in
Section~\ref{Sextinction} and then survival (resp., the nontriviality
of the upper invariant law) for sufficiently large $\lambda$ in
Section~\ref{Ssurvival}.

%s2.1 #&#
\subsection{Extinction}\label{Sextinction}

We prove lower bounds on $\lambda_{\mathrm{c}}$ and $\lambda'_{\mathrm{c}}$ in
the present
subsection and upper bounds in the next. We start with $\lambda_{\mathrm{c}}$.

%le5 #&#
\begin{lemma}[(Triviality of the upper invariant law)]\label{Lnutriv}
For $\lambda\leq1$, the upper invariant law of the cooperative
branching-coalescent satisfies $\bar\nu=\delta_\varnothing$.
\end{lemma}

\begin{pf}
In this proof, it will be more convenient to
work with the process
$(X_t)_{t\geq0}$ taking values in $\{0,1\}^{\mathbb Z}$, rather than using
set notation as in (\ref{etadef}).

Let $X$ be a cooperative branching-coalescent started in any
translation-invariant initial law. For any $x_0,\ldots,x_n\in\{0,1\}$,
let us write for $t\geq0$
%
%e12 #&#
\begin{equation}
\label{ptdef}\qquad p_t(x_0x_1\cdots
x_n):={\mathbb P} \bigl[X_t(i)=x_0,X_t(i+1)=x_1,
\ldots,X_t(i+n)=x_n \bigr],
\end{equation}
which does not depend on $i\in{\mathbb Z}$ by the translation
invariance of our process
and the initial law. It follows from basic generator calculations that
%
%e13 #&#
\begin{eqnarray}\label{difp}
\frac{\partial}{\partial{t}}p_t(1)&=&-p_t(1)+
\frac{1}{2}p_t(10)+\frac{1}{2}p_t(01) +
\frac{1}{2}\lambda p_t(110)+\frac{1}{2}\lambda
p_t(011)\nonumber
\\
&=&-p_t(11)+\lambda\bigl(p_t(11)-p_t(111)
\bigr)
\\
&=&(\lambda-1)p_t(11)-\lambda p_t(111).\nonumber
\end{eqnarray}
Here, the terms in the first line arise from a particle at $i$ jumping
away as
well as a vacant site at $i$ becoming occupied by particles jumping
there or
by pairs of particles giving birth to a particle at $i$. We have rewritten
this using that $p_t(1)=p_t(10)+p_t(11)$ and
$p_t(11)=p_t(110)+p_t(111)$, and
similar relations for $p_t(01)$ and $p_t(011)$.

Now imagine that $X_0$ is distributed according to $\bar\nu$, or
in fact any
translation invariant stationary law. Then, assuming moreover that
$0<\lambda\leq
1$, we have
%
%e14 #&#
\begin{equation}
0=\frac{\partial}{\partial{t}}p_t(1)\leq-\lambda p_t(111)=-\lambda
p_0(111),
\end{equation}
from which we conclude that
%
%e15 #&#
\begin{equation}
\label{p111} p_0(111)={\mathbb P} \bigl[ \bigl(X_0(1),X_0(2),X_0(3)
\bigr)=(1,1,1) \bigr]=0.
\end{equation}
We will show that this implies that $X_0$ is identically zero a.s.

Indeed, if $X_0$ is not identically zero, then by translation invariance
$p_t(1)=p_0(1)={\mathbb P}[X_0(i)=1]=:\varepsilon>0$ $(i\in{\mathbb
Z})$ so for $n>3\varepsilon^{-1}$ the
expected number of particles in $\{1,\ldots,n\}$ is greater than
three. In
particular, there is a positive probability of finding three particles
in this
interval. Using Lemma~\ref{Lmonot}, we may estimate $X$ from below by a
system of coalescing random walks without cooperative branching. Since there
is a positive probability that three coalescing random walks started anywhere
in $\{1,\ldots,n\}$ end up at the sites $1,2,3$ at time 1, using stationarity
we see that the probability in (\ref{p111}) is positive, contradicting our
assumption.

If $\lambda=0$, then the same argument applies, except that we use that
$0=\frac{\partial}{\partial{t}}p_t(1)=p_t(11)$, and we only need to
show that this implies the
triviality of $X$, which is weaker than what we have already shown. Our
arguments show that for $\lambda\leq1$, no translation invariant
stationary law
can exist that is not concentrated on the empty configuration. In particular,
the upper invariant law must be concentrated on the empty
configuration.
\end{pf}

%le6 #&#
\begin{lemma}[(Extinction)]\label{Lextinct}
For $\lambda\leq1$, the cooperative branching-coalescent started in
any finite,
nonempty initial state $A$ satisfies
%
%e16 #&#
\begin{equation}
{\mathbb P} \bigl[\exists T<\infty\mbox{ s.t. }\bigl|\eta^A_t\bigr|=1
\ \forall t\geq T \bigr]=1.
\end{equation}
\end{lemma}

\begin{pf}
Given $\eta^A_t$ we have that $|\eta^A_t|$
increases by $1$ due to cooperative branching at rate
%
%e17 #&#
\begin{equation}
\frac{\lambda}{2} \sum_{i \in{\mathbb Z}} ( 1_{\{ \{i, i+1\}
\subset\eta^A_t, i+2 \notin\eta^A_t\}}
+1_{\{ \{i, i+1\} \subset\eta^A_t, i-1 \notin\eta^A_t\}} )
\end{equation}
and decreases by $1$ due to coalescence at rate $\sum_{i \in{\mathbb
Z}} 1_{\{ \{i, i+1\} \subset\eta^A_t\}}$.
Therefore, we obtain
%
%e18 #&#
\begin{eqnarray}
\frac{\partial}{\partial{t}}{\mathbb E} \bigl[\bigl|\eta^A_t\bigr| \bigr] &=&
\frac{\lambda}{2} \sum_{i\in{\mathbb Z}} \bigl( {\mathbb P}\bigl[
\{ i,i+1\}\subset\eta^A_t, i+2 \notin\eta^A_t
\bigr]\nonumber
\\
&&\hspace*{26pt}{}
 +{\mathbb P}\bigl[\{i,i+1\}\subset\eta^A_t,
i-1 \notin\eta^A_t \bigr] \bigr)
\\
& &{}-\sum_{i\in{\mathbb Z}}{\mathbb P}\bigl[\{i,i+1\}
\subset\eta^A_t \bigr].\nonumber
\end{eqnarray}
Since
%
%e19 #&#
\begin{equation}
1_{\{ \{i, i+1\} \subset\eta^A_t, i+2 \notin\eta^A_t\}}= 1_{\{ \{i,
i+1\} \subset\eta^A_t\}} -1_{\{ \{i, i+1, i+2\} \subset\eta^A_t\}},
\end{equation}
it follows from a calculation as in (\ref{difp}) using the translation
invariance that
%
%e20 #&#
\begin{eqnarray}
\frac{\partial}{\partial{t}}{\mathbb E} \bigl[\bigl|\eta^A_t\bigr| \bigr] &=& (
\lambda-1)\sum_{i\in{\mathbb Z}}{\mathbb P}\bigl[\{i,i+1\}\subset
\eta^A_t \bigr]
\nonumber\\[-8pt]\\[-8pt]
&&{} -\lambda\sum
_{i\in{\mathbb Z}}{\mathbb P}\bigl[\{i,i+1,i+2\}\subset\eta
^A_t \bigr].\nonumber
\end{eqnarray}
In particular, if $\lambda\leq1$ this is easily seen to imply due to
the Markov property that
$|\eta^A_t|$ is a supermartingale with respect to $ {\mathcal
F}_t^A:=\sigma( \eta^A_s, 0 \leq s \leq t)$ since for $0 \leq s \leq t$,
%
%e21 #&#
\begin{equation}
\qquad {\mathbb E} \bigl[\bigl|\eta^A_t\bigr| | {\mathcal
F}_s^A \bigr] = {\mathbb E} \bigl[\bigl|\eta^{\eta^A_s}_{t-s}\bigr|
| \eta_s^A \bigr] = \bigl| \eta_s^A\bigr|
+ \int_0^{t-s} \frac{\partial}{\partial{u}}{\mathbb E}
\bigl[\bigl|\eta^{ \eta_s^A }_u\bigr| | \eta_s^A
\bigr] \,\mathrm{d}u \leq\bigl| \eta_s^A\bigr|.
\end{equation}
By supermartingale convergence, it
follows that
%
%e22 #&#
\begin{equation}
\label{ascon} \bigl|\eta^A_t\bigr|\mathop{\longrightarrow}_{{t}\to\infty}
N\qquad\mbox{a.s.}
\end{equation}
for some ${\mathbb N}$-valued random variable $N$. Let
%
%e23 #&#
\begin{equation}
{\mathcal A}_T:= \bigl\{\exists t\geq T\mbox{ s.t. }\bigl|
\eta^A_{t-}\bigr|\neq\bigl|\eta^A_t\bigr| \bigr\}
\end{equation}
denote the event that the number of particles will change at some time greater
or equal than $T$, and let $\rho(A)$ denote the probability of
${\mathcal A}_0$ as a
function of the initial state $A$. Using the continuity of
conditional probabilities w.r.t. the $\sigma$-field (see \cite{Chu74},
Theorem~9.4.8,
or \cite{Bil86}, Theorems~3.5.5 and 3.5.7), we conclude that for each
$S\leq T$,
%
%e24 #&#
\begin{equation}
\rho\bigl(\eta^A_T\bigr)={\mathbb P}\bigl[{\mathcal
A}_T|{\mathcal F}^{A}_T\bigr]\leq{\mathbb P}
\bigl[{\mathcal A}_S|{\mathcal F}^{A}_T\bigr]
\mathop{\longrightarrow}_{{T}\to\infty} {\mathbb P}\bigl[{\mathcal
A}_S|{\mathcal F}^{A}_\infty\bigr]=1_{{\mathcal
A}_S}
\qquad\mbox{a.s.}
\end{equation}
It follows that $\lim_{T\to\infty}\rho(\eta^A_T)=0$ a.s. on the
complement of
the event\break $\bigcap_{S\geq0}{\mathcal A}_S$; that is, the event
%
%e25 #&#
\begin{equation}
\lim_{T\to\infty}\rho\bigl(\eta^A_T\bigr)=0
\quad\mbox{or}\quad\forall S\geq0\ \exists t\geq S\mbox{ s.t. } \bigl|
\eta^A_{t-}\bigr|\neq\bigl|\eta^A_t\bigr|
\end{equation}
has probability one. By (\ref{ascon}), we conclude that
$\lim_{T\to\infty}\rho(\eta^A_T)=0$ a.s. By the recurrence of
one-dimensional
random walk, it is easy to see that $\rho$ is uniformly bounded away
from zero
on $\{A\dvtx |A|\geq2\}$
(in fact, it is not hard to see that $\rho\equiv1$ on
this set), so we conclude that $\lim_{T\to\infty}|\eta^A_T|=1$ a.s.
\end{pf}

%s2.2 #&#
\subsection{Survival}\label{Ssurvival}

In this section we show that for $\lambda$ sufficiently large, the cooperative
branching-coalescent survives and has a nontrivial upper invariant law.
As a
first step, we compare it from below with a contact process with ``double
deaths.'' Since in the cooperative branching-coalescent, only pairs of
particles can produce offspring, we wish to estimate the number of occupied
neighboring pairs from below.

For\vspace*{2pt} each $i\in{\mathbb Z}+\frac{1}{2}$, let $\rvec\pi(i),\lvec\pi
(i),\pi^\ast(i)$ be
independent Poisson subsets of ${\mathbb R}$ with intensities
$\frac{{1}}{{2}}\lambda,\frac{{1}}{{2}}\lambda$, and $1$,
respectively. For each
$\zeta_0\subset{\mathbb Z}$, we may construct a Markov process
$(\zeta_t)_{t\geq0}$ with
initial state $\zeta_0$ that evolves according to the following rules.

For each $i\in{\mathbb Z}+\frac{1}{2}$, if immediately prior to some
cooperative branching event
$t\in\rvec\pi(i)$ the state is $\zeta_{t-}$ and $i-\frac{1}{2}\in
\zeta_{t-}$, then we
set $\zeta_t:=\zeta_{t-}\cup\{i+\frac{1}{2}\}$. If $i-\frac
{1}{2}\notin\zeta_{t-}$, then we do
nothing. A similar rule applies to $t\in\lvec\pi(i)$, where now the site
$i+\frac{1}{2}$, if occupied, infects the site $i-\frac{1}{2}$.
Finally, for each
$i\in{\mathbb Z}+\frac{1}{2}$, at each time $t\in\pi^\ast(i)$, we replace
$\zeta_{t-}$ by $\zeta_t:=\zeta_{t-}\setminus\{i-\frac
{1}{2},i+\frac{1}{2}\}$.

With these rules, we see that $(\zeta_t)_{t\geq0}$ is a contact
process with
``double deaths,'' where sites infect each of their neighbors with infection
rate $\frac{1}{2}\lambda$, and for each pair $\{i,i+1\}$ of
neighboring sites, any
particles located at these sites die simultaneously with rate $1$.

%le7 #&#
\begin{lemma}[(Comparison with contact process with double deaths)]\label{Lcdd}
Let $(\eta_t)_{t\geq0}$ be a cooperative branching-coalescent with
cooperative branching rate $\lambda$, and let $(\zeta_t)_{t\geq0}$
be a contact
process with double deaths and infection rate $\frac{1}{2}\lambda$. Let
%
%e26 #&#
\begin{equation}
\eta^{(2)}_t:= \bigl\{i\in{\mathbb Z}\dvtx \{i,i+1\}\subset
\eta_t \bigr\} \qquad(t\geq0)
\end{equation}
denote the set of locations where $\eta_t$ contains a pair of neighboring
particles. Then $(\eta_t)_{t\geq0}$ and $(\zeta_t)_{t\geq0}$ can be coupled
such that
%
%e27 #&#
\begin{equation}
\label{zetet} \zeta_0\subset\eta^{(2)}_0
\quad\mbox{implies}\quad\zeta_t\subset\eta^{(2)}_t
\qquad(t\geq0).
\end{equation}
\end{lemma}

\begin{pf}
We claim that (\ref{zetet}) holds if we
construct $(\eta_t)_{t\geq
0}$ by means of a graphical representation with Poisson sets
$\lvec\omega(i),\rvec\omega(i)$ as in Section~\ref{Sdefs} and construct
$(\zeta_t)_{t\geq0}$ by means of a graphical representation with
Poisson sets
given by
%
%e28 #&#
\begin{eqnarray}
\lvec\pi\bigl(i-\tfrac{1}{2}\bigr)&:=& \lvec\omega(i),\qquad\rvec\pi
\bigl(i-\tfrac
{1}{2}\bigr):=\rvec\omega(i),
\nonumber\\[-8pt]\\[-8pt]
\pi^\ast \bigl(i-\tfrac{1}{2}\bigr)&:=&\lvec\omega\bigl(i-\tfrac{1}{2}\bigr)
\cup\rvec\omega\bigl(i+\tfrac{1}{2}\bigr),\nonumber
\end{eqnarray}
$(i\in{\mathbb Z})$, which are independent Poisson sets with intensities
$\frac{{1}}{{2}}\lambda,\frac{{1}}{{2}}\lambda$ and $1$, respectively.

It suffices to check that if $\zeta_t\subset\eta^{(2)}_t$ is true
just prior to a
cooperative branching event or coalescing jump event, then it will also be
true immediately after such an event. For $i\in{\mathbb Z}$, if prior
to some
$t\in\rvec\omega(i)=\lvec\pi(i-\frac{1}{2})$ one has $\{i-1,i\}
\subset\eta_{t-}$ and
$i-1\in\zeta_{t-}$, then $\zeta_t=\zeta_{t-}\cup\{i\}$ while now also
$\{i,i+1\}\subset\eta_t$ since the pair $\{i-1,i\}$ has given birth
to a particle
at $i+1$. The same argument applies to cooperative branching events to the
left. For $i\in{\mathbb Z}$, it may happen that a pair $\{i,i+1\}
\subset\eta_{t-}$ is
destroyed due to a coalescing jump event
%
%e29 #&#
\begin{eqnarray}
t&\in&\bigl(\lvec\omega\bigl(i-\tfrac{1}{2}\bigr)\cup\rvec\omega\bigl(i+
\tfrac
{1}{2}\bigr) \bigr) \cup\bigl(\lvec\omega\bigl(i+\tfrac{1}{2}
\bigr)\cup\rvec\omega\bigl(i+\tfrac
{{3}}{{2}}\bigr) \bigr)
\nonumber\\[-8pt]\\[-8pt]
& =&
\pi^\ast\bigl(i-\tfrac{1}{2}\bigr)\cup\pi^\ast\bigl(i+
\tfrac{1}{2}\bigr),\nonumber
\end{eqnarray}
which corresponds to the particle at $i$ or $i+1$ jumping to the left or
right. But in this case, $i\notin\zeta_t$ since any particles
on either $\{i-1,i\}$ or $\{i,i+1\}$ have died simultaneously. Coalescing
jump events may also lead to the creation of new pairs but also in this case,
the inclusion $\eta^{(2)}_t\supset\zeta_t$ is preserved.
\end{pf}

Clearly, if the contact process with double deaths $(\zeta_t)_{t\geq
0}$ with
infection rate $\frac{{1}}{{2}}\lambda$ survives, then so does the cooperative
branching-coalescent with cooperative branching rate $\lambda$. We
note that
numerical simulations indicate that the contact process with double
deaths has
a critical infection rate of approximately $3.65\pm0.05$, so
presumably this
happens for approximately $\lambda\geq7.3\pm0.1$. The contact
process with double
deaths is a monotone particle system, so by the same arguments as for the
cooperative branching-coalescent [see (\ref{upper})], it has an upper
invariant law. Coupling the processes $\eta^{\mathbb Z}_t$ and $\zeta
^{\mathbb Z}_t$ as in
Lemma~\ref{Lcdd} and sending $t\to\infty$, we see that if the contact
process with double deaths has a nontrivial upper invariant law, then
so does
the cooperative branching-coalescent.

Thus we are left with the task of proving that for sufficiently large
infection rate $\frac{1}{2}\lambda$, the contact process with double deaths
survives and has a nontrivial upper invariant law. In fact, it suffices to
prove the first statement only. This is because the contact process with
double deaths is self-dual in the sense that ${\mathbb P}[\zeta
^A_t\cap
B\neq\infty]={\mathbb P}[A\cap\zeta^B_t\neq\infty]$ $(A,B\in
{\mathbb Z}, t\geq0)$, just like
the normal contact process (as can easily be proved from the graphical
representation), and hence its upper invariant law is nontrivial (for a given
value of $\lambda$) if and only if the process survives. (See the
discussion for
the standard contact process around formulas (I.1.7) and (I.1.8) in
\cite{Lig99}.)

Unfortunately, there seems to be no easy way to compare the contact process
with double deaths with a normal contact process. There exist several
ways of
proving survival (for sufficiently large $\lambda$) of the standard,
one-dimensional contact process. Each of these might be attempted for the
contact process with double deaths as well. We will use the most robust
technique, comparison with oriented percolation, which, however,
performs rather
poorly when it comes to finding explicit upper bounds on the critical
value. We will not attempt to find such explicit bounds.

Let\vspace*{1pt} ${\mathbb Z}^2_{\mathrm{even}}:=\{(i,n)\in{\mathbb Z}^2\dvtx i+n\mbox{ is
even}\}$. We equip
${\mathbb Z}^2_{\mathrm{even}}$ with the structure of a directed graph by
drawing for each
$z=(i,n)\in{\mathbb Z}^2_{\mathrm{even}}$ two directed edges (arrows)
$e^-_z$ and $e^+_z$
which point from $(i,n)$ to $(i-1,n+1)$ and $(i+1,n+1)$, respectively. Let
$E$ be the set of all directed edges $e^\pm_z$, and let $(\chi
_e)_{e\in E}$ be
i.i.d. Bernoulli random variables with ${\mathbb P}[\chi_e=1]=p$. We
say that the
edge $e$ is open if $\chi_e=1$. For $z,z'\in{\mathbb Z}^2_{\mathrm{even}}$, we say that
there is an open path from $z$ to $z'$, denoted as $z\to z'$, if either $z=z'$
or $z=(i,n)$, $z'=(i',n')$ with $n'>n$ and there exists a function
$\gamma\dvtx \{n,\ldots,n'\}\to{\mathbb Z}$ such that $\gamma_n=i$,
$\gamma_{n'}=i'$, and for all
$k=n+1,\ldots,n'$ one has $|\gamma_k-\gamma_{k-1}|=1$ and the edge from
$(\gamma_{k-1},k-1)$ to $(\gamma_k,k)$ is open. For given $W_0\subset
{\mathbb Z}_{\mathrm{even}}$, we put
for $n\geq1$,
%
%e30 #&#
\begin{equation}
\label{W} W_n:=\bigl\{i\in{\mathbb Z}\dvtx (i,n)\in{\mathbb
Z}^2_{\mathrm{even}}, \exists i'\in W_0
\mbox{ s.t. }\bigl(i',0\bigr)\to(i,n)\bigr\}.
\end{equation}
Then $W=(W_n)_{n\geq0}$ is a Markov chain, taking values, in turn, in the
subsets of ${\mathbb Z}_{\mathrm{even}}$ and ${\mathbb Z}_{\mathrm{odd}}$. We
call $W$ the \emph{oriented
percolation process}.

%pr8 #&#
\begin{proposition}[(Comparison with oriented percolation)]\label{Torcomp}
Let $(\zeta_t)_{t\geq0}$ denote the contact process with double
deaths, and
let $(W_n)_{n\geq0}$ denote the oriented percolation process with parameter
$p$. Then, for each $p<1$, there exists $\lambda',T>0$ such that for all
$\lambda\geq\lambda'$, the process $(\zeta_t)_{t\geq0}$ with
infection rate
$\frac{{1}}{{2}}\lambda$ and $(W_n)_{n\geq0}$ with parameter $p$ can
be coupled in
such a way that
%
%e31 #&#
\begin{equation}
\label{percomp} W_0\subset\zeta_0 \quad\mbox{implies}\quad
W_n\subset\zeta_{nT}\qquad(n\geq0).
\end{equation}
\end{proposition}

\begin{pf}
We construct $(\zeta_t)_{t\geq0}$ using its
graphical
representation. By a trivial rescaling of time, we may assume that infection
events, corresponding to the Poisson sets $\lvec\pi(i),\rvec\pi
(i)$, have
intensity $\frac{1}{2}$ each, while deaths, corresponding to the
Poisson sets
$\pi^\ast(i)$, have intensity $\lambda^{-1}$. For each $T>0$, we
define a
collection of Bernoulli random variables $(\chi^T_e)_{e\in E}$ indexed
by the
edges of the directed graph $({\mathbb Z}^2_{\mathrm{even}},E)$, in the
following way. For
the directed edge $e^+_z$ from $z=(i,n)$ to $(i+1,n+1)$, we let $\chi
^T_e$ be
the indicator of the event
%
%e32 #&#
\begin{eqnarray}
\label{goodevent} & & \bigl\{\rvec\pi\bigl(i+\tfrac{1}{2}\bigr)
\cap\bigl(nT,(n+1)T\bigr]\neq\varnothing,
\nonumber\\[-8pt]\\[-8pt]
&&\qquad\bigl(\pi^\ast\bigl(i-\tfrac{1}{2}
\bigr)\cup\pi^\ast\bigl(i+\tfrac{1}{2}\bigr)\cup\pi
^\ast\bigl(i+\tfrac{{3}}{{2}}\bigr) \bigr) \cap\bigl(nT,(n+1)T\bigr]=
\varnothing\bigr\},\nonumber
\end{eqnarray}
which says that there is an infection from $i$ to $i+1$ in the time interval
$(nT,(n+1)T]$, but no deaths occur in $i$ or $i+1$ during this time interval.
For directed edges $e^-_z$ to the left, the analogous definition
applies. Clearly, if ${\mathbb Z}_{\mathrm{even}}\ni i\in\zeta_0$ and
there exists a path
from $(i,0)$ to $(i',n)$ along edges that are open in the sense of the
$(\chi^T_e)_{e\in E}$, then $i'\in\zeta_{nT}$.

By first choosing $T$ large enough and then $\lambda$ large enough we
can make the
probability of the event in (\ref{goodevent}) as close to one as we
wish. The
events belonging to different edges are not independent, but they are
$m$-dependent for a suitable $m$, so by standard results
\cite{Lig99}, Theorem~B26, the $(\chi^T_e)_{e\in E}$ can be estimated
from below
by i.i.d. Bernoulli random variables with a succes probability $p$
that can
be made arbitrarily close to one. Using these i.i.d. Bernoulli random
variables to construct the oriented percolation process, we arrive at
(\ref{percomp}).
\end{pf}

\begin{pf*}{Proof of Theorem~\ref{Tphase}} The facts that
the upper invariant law is
trivial and the process dies out for $\lambda\leq1$ have been proved in
Lemmas~\ref{Lnutriv} and \ref{Lextinct}. To prove that for $\lambda
$ suffiently
large, the upper invariant law is nontrivial and the process survives, by
Lemma~\ref{Lcdd} and the discussion below it, it suffices to show
that the
contact process with double deaths survives for $\lambda$ suffiently
large. This
follows from Theorem~\ref{Torcomp} and the fact that the oriented percolation
process $(W_n)_{n\geq0}$ survives for $p>8/9$ by~\cite{Dur88},
Section~5a.
\end{pf*}

%s3 #&#
\section{Decay of the density}\label{Sdecay}

%s3.1 #&#
\subsection{Some general terminology}\label{Stermin}

Recall (see, e.g., \cite{Lig85}, Definition~II.3.1) that two Markov
processes $X$ and $Y$ with metrizable state spaces
$S_X$ and $S_Y$ are dual to each other with bounded, Borel measurable duality
function $\psi\dvtx S_X\times S_Y\to{\mathbb R}$, if for processes with arbitrary
deterministic initial states $X_0$ and $Y_0$ one has
%
%e33 #&#
\begin{equation}
\label{dual} {\mathbb E}\bigl[\psi(X_t,Y_0)\bigr]={
\mathbb E}\bigl[\psi(X_0,Y_t)\bigr]\qquad(t\geq0).
\end{equation}
If (\ref{dual}) holds for deterministic initial states, then it holds more
generally when $X$ and $Y$ are independent and have (possibly) random initial
states, as can be seen by integrating both sides of (\ref{dual})
w.r.t. the
product of the laws of $X_0$ and $Y_0$.

More generally, borrowing terminology from \cite{AS05}, we say that
$Y$ is a
\emph{subdual} of $X$ if
%
%e34 #&#
\begin{equation}
\label{subdual} {\mathbb E}\bigl[\psi(X_t,Y_0)\bigr]\geq{
\mathbb E}\bigl[\psi(X_0,Y_t)\bigr]\qquad(t\geq0)
\end{equation}
whenever $X$ and $Y$ are independent. In particular, if $Y$ is started
in an
invariant law, and hence $(Y_t)_{t\geq0}$ is a stationary process,
then this
implies that the function
%
%e35 #&#
\begin{equation}
h(x):={\mathbb E}\bigl[\psi(x,Y_t)\bigr]\qquad(x\in
S_X, t\geq0)
\end{equation}
is a subharmonic function for the process $X$. We define superduals (which
then may give rise to superharmonic functions) similarly, by reversing the
inequality sign.

Following \cite{JK12}, we say that a duality as in (\ref{dual}) is a
\emph{pathwise} duality if for each $t>0$, it is possible to couple the processes
$X$ and $Y$, which have c\`adl\`ag sample paths, in such a way that the
stochastic
process
%
%e36 #&#
\begin{equation}
\label{pathwise} s\mapsto\psi(X_{s-},Y_{t-s})
\end{equation}
is a.s. constant on $[0,t]$. Likewise, we may say that we have a pathwise
subduality (resp., superduality) if this function is a.s. nondecreasing
(nonincreasing).

%s3.2 #&#
\subsection{Coalescing random walk duality}\label{Scodu}

In this section, we consider the case that the cooperative branching
rate is
zero. In this case, the cooperative branching-coalescent $\eta^A_t$ from
(\ref{etadef}) reduces to a system of coalescing random walks, and the
graphical representation contains only coalescing jump events, that is,
$\lvec\omega(i)=\varnothing=\rvec\omega(i)$ for each $i\in{\mathbb Z}$.

By definition, an \emph{open path} in our graphical representation is
a c\`adl\`ag
function $\xi\dvtx L\to{\mathbb Z}$, defined on some interval $L\subset
{\mathbb R}$, satisfying the
following rules:
\begin{longlist}[(2)]
\item[(1)] If $t\in\lvec\omega(\xi_{t-}-\frac{1}{2})$ [resp., $t\in\rvec
\omega(\xi_{t-}+\frac{1}{2})$]
for some $t\in L$, then $\xi_t=\xi_{t-}-1$ (resp., $=\xi_{t-}+1$).
\item[(2)] If for some $t\in L$,
$t\notin(\lvec\omega(\xi_{t-}-\frac{1}{2})\cup\rvec\omega(\xi
_{t-}+\frac{1}{2}))$, then
$\xi_t=\xi_{t-}$.
\end{longlist}

In words, this says that $\xi$ walks upwards until it meets the rear
end of an
arrow, at which instance it jumps to the tip of the arrow and continues its
journey upwards. For each deterministic \emph{starting point}
$(i,s)\in{\mathbb Z}\times{\mathbb R}$, there a.s. exists a unique
open path
$\xi^{(i,s)}\dvtx [s,\infty)\to{\mathbb Z}$ such that $\xi^{(i,s)}_s=i$,
and this path is
distributed as a random walk that jumps to the
positions immediately on its left or\vspace*{1pt} right with rate $\frac{1}{2}$
each. Moreover,
paths started from different starting points are independent until the first
time they meet and coalesce (i.e., go on as a single walker as soon as they
meet).

%
%f4 #&#
\begin{figure}[t]

\includegraphics{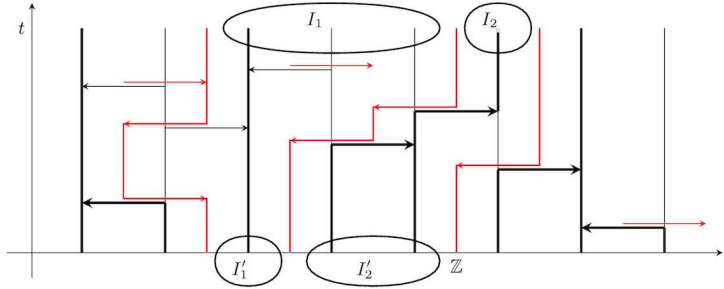}

%
%%need to make coordinates y1--y26
%%below 0.4 should be \y1
%
% -- (0.5,\vvw+0.05) -- (1.5,\vvw+0.05) -- (1.5,5.3);
% -- (3.5,\vvz+0.05)-- (2.5,\vvz+0.05)-- (2.5,5.3);
% -- (4.5,5.3);
%
%
%coordinates{(1.7,8) (3,8.3) (4.3,8) (3,7.7)};
%coordinates{(4.6,8) (5,8.3) (5.4,8) (5,7.7)};
%
%coordinates{(1.6,5.2) (2,5.5) (2.4,5.2) (2,4.9)};
%coordinates{(2.7,5.2) (3.5,5.5) (4.3,5.2) (3.5,4.9)};
\caption{Coalescing random walk duality. One has $\eta_t\cap I_1\neq
\varnothing$
if and only if $\eta_0\cap I'_1\neq\varnothing$, and $\eta_t\cap
I_2\neq\varnothing$
if and only if $\eta_0\cap I'_2\neq\varnothing$. The\vspace*{1pt} boundaries
of the intervals $I_1$ and $I_2$, respectively, $I_1'$~and~$I_2'$, are dual
coalescing random walk paths evaluated at time $t$, respectively,~$0$.}
\label{figcodu}
\end{figure}

It is well known (see, e.g., the Appendix of \cite{STW00} for an analogous
duality in discrete time)
that such systems of coalescing random walks
are self-dual in the following sense. Set
%
%e37 #&#
\begin{equation}
\hat{\lvec\omega}(i):=\rvec\omega\bigl(i+\tfrac{1}{2}\bigr)\quad\mbox
{and}\quad\hat{\rvec\omega}(i):=\lvec\omega\bigl(i-\tfrac{1}{2}\bigr)
\qquad(i\in{\mathbb Z}),
\end{equation}
and for each $t\in\hat{\lvec\omega}(i)$ [resp., $t\in\hat{\rvec
\omega}(i)$], draw a
\emph{dual arrow} from $i+\frac{1}{2}$ to $i-\frac{1}{2}$ (resp., from
$i-\frac{1}{2}$ to $i+\frac{1}{2}$). In
Figure~\ref{figcodu}, these dual arrows have been drawn in red.
By definition, a \emph{dual open path} in our graphical representation
is a
caglad (i.e., left continuous with right limits) function
$\hat\xi\dvtx L\to{\mathbb Z}+\frac{1}{2}$, defined on some interval
$L\subset{\mathbb R}$,
such that:
\begin{longlist}[(2)]
\item[(1)] If $t\in\hat{\lvec\omega}(\hat\xi_{t+}-\frac{1}{2})$
[resp., $t\in\hat{\rvec\omega}(\hat\xi_{t+}+\frac{1}{2})$] for
some $t\in L$, then
$\hat\xi_t=\hat\xi_{t+}-1$ (resp., $=\hat\xi_{t+}+1$).
\item[(2)] If for some $t\in L$,
$t\notin(\hat{\lvec\omega}(\hat\xi_{t+}-\frac{1}{2})\cup\hat
{\rvec\omega}(\hat\xi_{t+}+\frac{1}{2}))$,
then $\hat\xi_t=\hat\xi_{t+}$.
\end{longlist}
In words, this says that the dual open paths walk downwards in time
until they
meet the rear end of a dual arrow, at which instance they jump to its
tip and
continue their journey downwards. For each deterministic starting point
$(i,s)\in({\mathbb Z}+\frac{1}{2})\times{\mathbb R}$, there a.s. exists
a unique dual open path
$\hat\xi^{(i,s)}\dvtx (-\infty,s]\to{\mathbb Z}+\frac{1}{2}$ such that
\mbox{$\hat\xi^{(i,s)}_s=i$}, and these
downward paths are distributed in the same way as the forward paths,
except for
a rotation over 180 degrees.

We observe that forward and dual open paths do not cross. As a result, for
each deterministic $t>0$ and $i,j\in{\mathbb Z}+\frac{1}{2}$ with
$i<j$, if there exists a
forward open path starting at time $0$ that passes at time $t$ between
$i$ and
$j$, then such a path must start between $\hat\xi^{(i,t)}_0$ and
$\hat\xi^{(j,t)}_0$; in particular, such a forward path can exist
only if
$\hat\xi^{(i,t)}_0<\hat\xi^{(j,t)}_0$. Conversely, any forward path that
starts between $\hat\xi^{(i,t)}_0$ and $\hat\xi^{(j,t)}_0$ must
pass at time
$t$ between $i$ and $j$. For any $i,j\in{\mathbb Z}+\frac{1}{2}$ with
$i<j$, let us write
%
%e38 #&#
\begin{equation}
\label{haint} \langle i,j\rangle:=\{k\in{\mathbb Z}\dvtx i<k<j\}=\bigl\{i+
\tfrac{1}{2},\ldots,j-\tfrac{1}{2}\bigr\}.
\end{equation}
Then, if $\eta^A$ is our system of coalescing random walks defined as in
(\ref{etadef}) (with cooperative branching rate $\lambda=0$), then by the
arguments we have just given,
%
%e39 #&#
\begin{equation}
\label{walkdual} \eta^A_t\cap\langle i,j\rangle\neq
\varnothing\quad\mbox{if and only if}\quad A\cap\bigl\langle\hat\xi^{(i,t)}_0,
\hat\xi^{(j,t)}_0\bigr\rangle\neq\varnothing\qquad\mbox{a.s.}
\end{equation}
(See Figure~\ref{figcodu}.)
In fact, the process
%
%e40 #&#
\begin{equation}
\bigl(\bigl\langle\hat\xi^{(i,t)}_{t-s},\hat
\xi^{(j,t)}_{t-s}\bigr\rangle\bigr)_{s\geq0}
\end{equation}
defines a nearest-neighbor voter model by specifying the clusters that
are occupied
by either 0's or 1's. Relationship (\ref{walkdual}) is then a special
case of
the well-known (pathwise) duality between coalescing random walks and voter
models; see, for example, \cite{Lig85}, Section~V.1.

For each $i,j\in{\mathbb Z}$, let
%
%e41 #&#
\begin{equation}
\label{tauij} \tau_{i,j}:=\inf\bigl\{t\geq0\dvtx \xi^{(i,0)}_t=
\xi^{(j,0)}_t\bigr\}
\end{equation}
denote the time at which $\xi^{(i,0)}$ and $\xi^{(j,0)}$ coalesce.
We will be especially interested in
%
%e42 #&#
\begin{equation}
\label{tau23} \tau_{\langle2\rangle}:=\tau_{0,1}\quad\mbox{and}\quad\tau
_{\langle3\rangle}:=\tau_{0,1}\wedge\tau_{1,2},
\end{equation}
which are the first time that any pair out of two (resp., three) walkers
meet each other and coalesce when the walkers are initially located at
neighboring positions.

Consider the system $\eta^{\mathbb Z}$ of coalescing random walks
started with each
site occupied. As in (\ref{ptdef}), let
%
%e43 #&#
\begin{equation}
p_t(1)={\mathbb P} \bigl[i\in\eta^{\mathbb Z}_t
\bigr]\quad\mbox{and}\quad p_t(11)={\mathbb P} \bigl[i\in
\eta^{\mathbb Z}_t, i+1\in\eta^{\mathbb Z}_t \bigr]
\end{equation}
denote the density of occupied sites and the density of pairs of occupied
neighboring sites, respectively, as a function of time. We claim that
%
%e44 #&#
\begin{equation}
\label{ptau} p_t(1)={\mathbb P}[t<\tau_{\langle2\rangle}]\quad
\mbox{and}\quad p_t(11)={\mathbb P}[t<\tau_{\langle3\rangle}].
\end{equation}
Indeed, by the coalescing random walk duality (\ref{walkdual}),
%
%e45 #&#
\begin{equation}
p_t(1)={\mathbb P} \bigl[i\in\eta^{\mathbb Z}_t
\bigr] ={\mathbb P} \bigl[\hat\xi^{(i-(1/2),t)}_0<\hat
\xi^{(i+(1/2),t)}_0 \bigr] ={\mathbb P}[t<\tau_{\langle2\rangle}],
\end{equation}
and similarly,
%
%e46 #&#
\begin{eqnarray}
p_t(11) &=& {\mathbb P} \bigl[i\in\eta^{\mathbb Z}_t,
i+1\in\eta^{\mathbb Z}_t \bigr]
\nonumber
\\[-8pt]
\\[-8pt]
&=&{\mathbb P} \bigl[\hat\xi^{(i-(1/2),t)}_0<\hat
\xi^{(i+(1/2),t)}_0<\hat\xi^{(i+(3/2),t)}_0 \bigr] ={
\mathbb P}[t<\tau_{\langle3\rangle}].
\nonumber
\end{eqnarray}

%s3.3 #&#
\subsection{Asymptotics of meeting times}

For any two functions $f,g\dvtx {[0,\infty)}\to(0,\infty)$, we write
%
%e47 #&#
\begin{equation}
f(t)\sim g(t)\qquad\mbox{as }t\to\infty
\end{equation}
to express the fact that
%
%e48 #&#
\begin{equation}
\frac{f(t)}{g(t)}\mathop{\longrightarrow}_{{t}\to\infty}1.
\end{equation}
Recall the definitions of $\tau_{\langle2\rangle}$ and $\tau
_{\langle3\rangle}$ from
(\ref{tau23}). We will need the following fact.

%le9 #&#
\begin{lemma}[(Asymptotics of meeting times)]\label{Lasmeet}
One has
%
%e49 #&#
\begin{equation}
\label{asmeet} \qquad{\mathbb P}[t<\tau_{\langle2\rangle}]\sim\frac{1}{\sqrt
{\pi}}t^{-1/2}
\quad\mbox{and}\quad{\mathbb P}[t<\tau_{\langle3\rangle}] \sim
\frac{1}{2\sqrt{\pi}}t^{-3/2} \qquad\mbox{as }t\to\infty.
\end{equation}
\end{lemma}

Note that the first statement about $\tau_{\langle2\rangle}$ is just
a result on the hitting time of zero of a random walk that is given
by the mutual distance of the two random walkers. The second statement
about $\tau_{\langle3\rangle}$ is less standard. For a proof
of both asymptotics in the case of a discrete time random walk, see,
for example, \cite{EK08}, Theorem 1.1. For completeness
we provide here a short proof of Lemma~\ref{Lasmeet}. For this we
need one preparatory technical lemma.

%le10 #&#
\begin{lemma}[(Asymptotic derivative)]\label{Lasder}
Let $\alpha>0$, let $F\dvtx {[0,\infty)}\to(0,\infty)$ be continuously
differentiable,
and assume that $t\mapsto\frac{\partial}{\partial{t}}F(t)$ is
nondecreasing. Assume moreover that
%
%e50 #&#
\begin{equation}
F(t)\sim t^{-\alpha}\qquad\mbox{as }t\to\infty.
\end{equation}
Then
%
%e51 #&#
\begin{equation}
\label{fsim} -\frac{\partial}{\partial{t}}F(t)\sim\alpha t^{-\alpha
-1}\qquad\mbox{as
}t\to\infty.
\end{equation}
\end{lemma}

\begin{pf}
Heuristically, we have
%
%e52 #&#
\begin{equation}
-\frac{\partial}{\partial{t}}F(t)\approx-\frac{\partial}{\partial
{t}}t^{-\alpha}=\alpha
t^{-\alpha-1},
\end{equation}
where it is not clear, a priori, how to interpret the approximate equality
$\approx$. It is easy to find examples showing that this cannot, in
general, be
interpreted in the sense of~$\sim$ without imposing further regularity
conditions (such as the monotonicity of the derivative).

To prove (\ref{fsim}), set $f(t):=-\frac{\partial}{\partial
{t}}F(t)$. We observe that for each
$\delta>0$,
%
%e53 #&#
\begin{eqnarray}
&&t^\alpha\int_t^{t(1+\delta)}f(s)\,
\mathrm{d}s\nonumber
\nonumber\\[-8pt]\\[-8pt]\nonumber
&&\qquad =t^\alpha F(t) -(1+\delta)^{-\alpha}\bigl(t(1+\delta)
\bigr)^\alpha F \bigl(t(1+\delta) \bigr)
\mathop{\longrightarrow}_{{t}\to\infty}1-(1+\delta)^{-\alpha}.\nonumber
\end{eqnarray}
Since $f$ is nonincreasing, it follows that
%
%e54 #&#
\begin{equation}
\liminf_{t\to\infty}(\delta t)t^\alpha f(t)\geq1-(1+\delta
)^{-\alpha}\qquad(\delta>0),
\end{equation}
and hence
%
%e55 #&#
\begin{equation}
\liminf_{t\to\infty}t^{\alpha+1}f(t)\geq\delta^{-1}
\bigl(1-(1+\delta)^{-\alpha}\bigr) \mathop{\longrightarrow}_{{\delta
\downarrow0}}
\alpha.
\end{equation}
In a similar fashion, by looking at the intergral from $t(1-\delta)$
to $t$,
we obtain that
%
%e56 #&#
\begin{equation}
\limsup_{t\to\infty}t^{\alpha+1}f(t)\geq\delta^{-1}
\bigl((1-\delta)^{-\alpha}-1\bigr) \mathop{\longrightarrow}_{{\delta
\downarrow0}}
\alpha.
\end{equation}\upqed
\end{pf}

\begin{pf*}{Proof of Lemma~\ref{Lasmeet}}
Let $(\Delta
_t)_{t\geq0}$ be a
continuous-time random walk on ${\mathbb Z}$, started in $\Delta_0=1$,
that jumps one step
to the left or right with rate one each, and let
%
%e57 #&#
\begin{equation}
\tau:=\inf\{t\geq0\dvtx \Delta_t=0\}.
\end{equation}
Then the distance between the two walkers $\xi^{(1,0)}_t-\xi
^{(0,0)}_t$ as a
function of time has the same distribution as $\Delta_t$ stopped at
$\tau$; in
particular $\tau_{\langle2\rangle}$ is equally distributed with
$\tau$. It is a
simple consequence of the reflection principle (compare
\cite{LPW09}, formula~(2.21)) that
%
%e58 #&#
\begin{equation}
{\mathbb P}[\Delta_t>0]={\mathbb P}[\Delta_t<0]+{\mathbb
P}[t<\tau]\qquad(t\geq0).
\end{equation}
Also, by symmetry, ${\mathbb P}[\Delta_t>2]={\mathbb P}[\Delta_t<0]$,
so we obtain that
%
%e59 #&#
\begin{equation}
{\mathbb P}[t<\tau]={\mathbb P} \bigl[\Delta_t\in\{1,2\} \bigr]\sim
\frac{1}{\sqrt{\pi}}t^{-1/2} \qquad\mbox{as }t\to\infty,
\end{equation}
where in the last step we have used the local central limit theorem
\cite{LL10}, Theorem~2.5.6, and the fact that $\var(\Delta_t)=2t$.

We recall from (\ref{ptau}) that ${\mathbb P}[t<\tau_{\langle
2\rangle}]$ and ${\mathbb P}[t<\tau_{\langle
3\rangle}]$ are given by the density of occupied sites $p_t(1)$ and the
density of pairs of occupied neighboring sites $p_t(11)$ in a system of
coalescing random walks started from the fully occupied state. By formula
(\ref{difp}) restricted to $\lambda=0$,
%
%e60 #&#
\begin{equation}
\frac{\partial}{\partial{t}}p_t(1)=-p_t(11). %\quad\mbox{and hence}\quad
\end{equation}
Applying Lemma~\ref{Lasder} we arrive at (\ref{asmeet}).
\end{pf*}

Lemma~\ref{Lasmeet} yields the following useful corollary.

%co11 #&#
\begin{corollary}[(Power-law bound)]\label{Cpbd}
\label{Lmeet2}
There exists a constant $K<\infty$ such that
%
%e61 #&#
\begin{equation}
\label{meet2} {\mathbb P}[t<\tau_{\langle3\rangle}]\leq Kt^{-3/2}\qquad(t
\geq0).
\end{equation}
\end{corollary}

%s3.4 #&#
\subsection{Mean meeting time of three walkers}

It follows from Lemma~\ref{Lasmeet} that ${\mathbb E}[\tau_{\langle
2\rangle}]=\infty$ but
${\mathbb E}[\tau_{\langle3\rangle}]<\infty$. In fact, it turns out
that the expectation of
$\tau_{\langle3\rangle}$ is exactly one. While this fact is not
essential in the
following, it simplifies our formulas and makes our estimates more explicit.
In view of this, we provide a proof here. Although the content of
Lemma~\ref{Lmeanmeet} below seems to be known, we did not find a reference.

Recall that for each $i\in{\mathbb Z}$, $(\xi^{(i,0)}_t)_{t\geq0}$
is a
continuous-time random walk on ${\mathbb Z}$ that jumps at the times of
a rate one
Poisson process to one of its neighboring sites, chosen with equal
probabilities. Walkers started at different sites jump independently
until they
meet, after which they coalesce. As in (\ref{tauij}), we let $\tau_{i,j}$
denote the first meeting time of the walkers started at $i$ and $j$.

%le12 #&#
\begin{lemma}[(Expected meeting time of three walkers)]\label{Lmeanmeet}
One has
%
%e62 #&#
\begin{equation}
\label{meanmeet} {\mathbb E}[\tau_{i,j}\wedge\tau_{j,k}]=(j-i)
(k-j)\qquad(i\leq j\leq k).
\end{equation}
\end{lemma}

\begin{pf}
Since we stop the process as soon as two
walkers meet, instead of looking at
coalescing random walks, we can equivalently study independent walkers. Let
$\vec\xi_t=(\xi^1_t,\xi^2_t,\xi^3_t)$ $(t\geq0)$ be three
independent walkers
started at $(\xi^1_0,\xi^2_0,\xi^3_0)=(i,j,k)$ with $i<j<k$. Then
$(\vec\xi_t)_{t\geq0}$ is a Markov process with generator
%
%e63 #&#
\begin{eqnarray}
Gf(i,j,k)&=&
\tfrac{{1}}{{2}} \bigl(f(i+1,j,k)+f(i-1,j,k)-2f(i,j,k) \bigr)\nonumber
\\
&&{}+\tfrac{{1}}{{2}} \bigl(f(i,j+1,k)+f(i,j-1,k)-2f(i,j,k)
\bigr)
\\
&&{}+\tfrac{{1}}{{2}} \bigl(f(i,j,k+1)+f(i,j,k-1)-2f(i,j,k)
\bigr).\nonumber
\end{eqnarray}
Set
%
%e64 #&#
\begin{eqnarray}
{\mathbb Z}^3_\leq &:=& \bigl\{(i,j,k)\in{\mathbb
Z}^3\dvtx i\leq j\leq k\bigr\},
\nonumber\\[-8pt]\\[-8pt]
{\mathbb Z}^3_<&:=&\bigl
\{(i,j,k)\in{\mathbb Z}^3\dvtx i<j<k\bigr\}.\nonumber
\end{eqnarray}
Consider the functions
%
%e65 #&#
\begin{equation}
\qquad f(i,j,k):=(j-i) (k-j)\quad\mbox{and}\quad h(i,j,k):=(j-i) (k-j) (k-i).
\end{equation}
Straightforward calculations give
%
%e66 #&#
\begin{equation}
\label{Gfh} Gf(i,j,k)=-1\quad\mbox{and}\quad Gh(i,j,k)=0\qquad\bigl((i,j,k)
\in{\mathbb Z}^3_<\bigr).
\end{equation}
%
%&=&\dis\ffrac{1}{2}\big(f(i-1,j,k)+f(i+1,j,k)+f(i,j-1,k)+f(i,j+1,k) %
%&=&\dis+f(i,j,k-1)+f(i,j,k+1)-6f(i,j,k)\big)=-1,
By Lemma~\ref{Lpolgro} below, the process
%
%e67 #&#
\begin{equation}
\label{MG} M_t:=f(\vec\xi_t)-\int
_0^t Gf(\vec\xi_s)\,\mathrm{d}s
\end{equation}
is a martingale with respect to the filtration generated by $\vec\xi$.
Therefore, setting
%
%e68 #&#
\begin{equation}
\tau:=\inf\bigl\{t\geq0\dvtx \vec\xi_t\notin{\mathbb Z}^3_<
\bigr\}=\tau_{i,j}\wedge\tau_{j,k}
\end{equation}
and using optional stopping, we see that for $\vec\xi_0=(i,j,k)\in
{\mathbb Z}^3_<$,
%
%e69 #&#
\begin{eqnarray}\label{ftau}
f(\vec \xi_0) &=& {\mathbb E}[M_{t\wedge\tau}] ={\mathbb E}
\bigl[f(\vec\xi_{t\wedge\tau})\bigr]-{\mathbb E}\biggl[\int_0^{t\wedge
\tau}Gf(
\vec\xi_s)\,\mathrm{d}s\biggr]
\nonumber\\[-8pt]\\[-8pt]
&=&{\mathbb E}\bigl[f(\vec\xi_{t\wedge\tau})\bigr]+{\mathbb E}[t
\wedge\tau],\nonumber
\end{eqnarray}
where we have used (\ref{Gfh}). Note that $\tau<\infty$ a.s. by the
recurrence of one-dimensional random walk. Therefore (\ref{meanmeet}) will
follow by letting $t\to\infty$ in (\ref{ftau}), provided we show that
%
%e70 #&#
\begin{equation}
\label{ftozero} {\mathbb E}\bigl[f(\vec\xi_{t\wedge\tau})\bigr]\mathop{
\longrightarrow}_{{t}\to\infty}0.
\end{equation}
Since $f$ is unbounded, this is not completely trivial. Indeed, our arguments
so far apply equally well to the function $f$ and the function $f':=f+h$,
while the right-hand side of (\ref{meanmeet}) is given by $f(i,j,k)$
and not
by $f'(i,j,k)$.
In order to prove~(\ref{ftozero}), we proceed as follows. Let
%
%e71 #&#
\begin{equation}
P_t(\vec\imath,\vec\jmath):={\mathbb P}^{\vec\imath}[\vec\xi
_{t\wedge\tau}=\vec\jmath]
\end{equation}
denote the transition probabilities of the stopped process.
Formula (\ref{Gfh}), Lemma~\ref{Lpolgro} below, and optional stopping
imply that
%
%e72 #&#
\begin{equation}
\bigl(h(\vec\xi_{t\wedge\tau})\bigr)_{t\geq0}
\end{equation}
is a martingale. As a result, setting
%
%e73 #&#
\begin{equation}
P^h_t(\vec\imath,\vec\jmath):=h(\vec
\imath)^{-1}P_t(\vec\imath,\vec\jmath)h(\vec\jmath) \qquad
\bigl(\vec\imath,\vec\jmath\in{\mathbb Z}^3_< \bigr)
\end{equation}
defines a transition probability on ${\mathbb Z}^3_<$. Let $\vec\xi
^h$ denote the
associated $h$-transformed Markov process, started in the same initial state
$\vec\xi^h_0=\vec\xi_0=(i,j,k)\in{\mathbb Z}^3_<$. Then, using the
fact that
$f=0$ on ${\mathbb Z}^3_\leq\setminus{\mathbb Z}^3_<$, we have that
%
%e74 #&#
\begin{eqnarray}
\label{hform} {\mathbb E}\bigl[f(\vec\xi_{t\wedge\tau})\bigr] &=& \sum
_{\vec\jmath\in{\mathbb Z}^3_\leq}P_t(\vec\imath,\vec\jmath)f(\vec
\jmath) =h(\vec\imath)\sum_{\vec\jmath\in{\mathbb Z}^3_<}P^h_t(
\vec\imath,\vec\jmath)f(\vec\jmath)h(\vec\jmath)^{-1}
\nonumber
\\[-8pt]
\\[-8pt]
&=&h(\vec\imath) {\mathbb E} \bigl[f\bigl(\vec\xi^h_t
\bigr)/h\bigl(\vec\xi^h_t\bigr) \bigr] =h(\vec\imath) {
\mathbb E} \bigl[\bigl(\xi^{h,3}_t-\xi^{h,1}_t
\bigr)^{-1} \bigr],
\nonumber
\end{eqnarray}
where we have used the notation
$\vec\xi^h_t=(\vec\xi^{h,1}_t,\ldots,\vec\xi^{h,3}_t)$.

We claim that $(\xi^{h,3}_t-\xi^{h,1}_t)^{-1}$ $(t\geq0)$ is a
supermartingale. Indeed, by optional stopping,\vspace*{2pt} the process
$(M_{t\wedge\tau})_{t\geq0}$ with $M$ as in (\ref{MG}) is a
martingale, so by
(\ref{Gfh}), $(f(\vec\xi_{t\wedge\tau}))_{t\geq0}$ is a supermartingale.
Setting $g(\vec\imath):=(i_3-i_1)^{-1}=f(\vec\imath)/h(\vec\imath
)$ $(\vec\imath\in{\mathbb Z}^3_<)$,
we see that
%
%e75 #&#
\begin{equation}
\qquad \sum_{\vec\jmath\in{\mathbb Z}^3_<}P^h_t(\vec
\imath,\vec\jmath)g(\vec\jmath) %=\sum_{\vecj\in\Z^3_<}h(\veci)^{-1}h(
%P^h_t(\veci,\vecj)h(\vecj)^{-1}f(\vecj)
=h(\vec
\imath)^{-1}\sum_{\vec\jmath\in{\mathbb Z}^3_<}P_t(\vec
\imath,\vec\jmath)f(\vec\jmath) \leq h(\vec\imath)^{-1}f(\vec\imath)=g(
\vec\imath).
\end{equation}
Since $(\xi^{h,3}_t-\xi^{h,1}_t)^{-1}$ is a\vspace*{2pt} bounded supermartingale, it
converges a.s. It is not hard to see that $\xi^{h,3}_t-\xi^{h,1}_t$ cannot
converge to a finite limit, so we conclude that
%
%e76 #&#
\begin{equation}
\bigl(\xi^{h,3}_t-\xi^{h,1}_t
\bigr)^{-1}\mathop{\longrightarrow}_{{t}\to
\infty}0\qquad\mbox{a.s.},
\end{equation}
which by (\ref{hform}) implies (\ref{ftozero}).
\end{pf}

\begin{rem*}
In a discrete-time setting, it is proved in
\cite{EK08}, Theorem 1.1(vi), that
%
%e77 #&#
\begin{equation}
t^{-{1}/{2}}\vec\xi^h_t\mathop{
\Longrightarrow}_{{t}\to\infty}B,
\end{equation}
with $B=(B^1,B^2,B^3)$ a random vector with density proportional to
%
%e78 #&#
\begin{equation}
\exp\bigl(-\tfrac{{1}}{{2}}\bigl(x_1^2 +
x_2^2 + x_3^2\bigr)\bigr)
h(x)^2.
\end{equation}
By looking at the associated jump chains, this result may be
transferred to
our present continuous-time setting.
\end{rem*}

We\vspace*{2pt} conclude this section by supplying the still outstanding lemma on
martingales. For each $\vec\imath=(i_1,\ldots,i_d)\in{\mathbb
Z}^d$, set
$\|\vec\imath\|:=\sup_{\alpha=1}^d|i_\alpha|$. We say that a
function $f\dvtx {\mathbb Z}^d\to{\mathbb R}$ is of
\emph{polynomial growth} if
%
%e79 #&#
\begin{equation}
\bigl|f(\vec\imath)\bigr|\leq K\bigl(1+\|\vec\imath\|^k\bigr)\qquad\bigl(\vec
\imath\in{\mathbb Z}^d\bigr)
\end{equation}
for some integers $K,k$.

%le13 #&#
\begin{lemma}[(Martingale problem for random walk)]\label{Lpolgro}
Let $(\vec\xi_t)_{t\geq0}$ be a conti\-nuous-time, nearest-neighbor
random walk
on ${\mathbb Z}^d$ started in a deterministic initial state, and let
$G$ denote its
generator. Then, for any function $f$ of polynomial growth, the process
%
%e80 #&#
\begin{equation}
M^f_t:=f(\vec\xi_t)-\int
_0^t Gf(\vec\xi_s)\,\mathrm{d}s
\qquad(t\geq0)
\end{equation}
is a martingale with respect to the filtration generated by
$(\vec\xi_t)_{t\geq0}$.
\end{lemma}

\begin{pf*}{Proof (sketch)}
Set $f_k(\vec\imath):=f(\vec
\imath)$ if $\|\vec\imath\|\leq k$ and $:=0$
otherwise. The fact that $M^{f_k}$ is a martingale is standard, so it suffices
to show that $M^{f_k}_t$ converges to $M^f_t$ in $L_1$-norm for each
$t\geq
0$. Now
%
%e81 #&#
\begin{equation}
\qquad {\mathbb E} \bigl[\bigl|M^f_t-M^{f_k}_t\bigr|
\bigr] \leq{\mathbb E} \bigl[\bigl|f(\vec\xi_t)-f_k(\vec
\xi_t)\bigr| \bigr] +\int_0^t{\mathbb E}
\bigl[\bigl|Gf(\vec\xi_s)-Gf_k(\vec\xi_s)\bigr|
\bigr]\,\mathrm{d}s.
\end{equation}
It is not hard to check that if $f$ is of polynomial growth, then so
are $Gf$
and $Gf_k$. Thus, $|f-f_k|$ and $|Gf-Gf_k|$ can be estimated by some function
of the form $K(1+\|\vec\imath\|^k)$, and the result follows by
dominated convergence
and the fact that nearest-neighbor random walk has moments of all orders.
\end{pf*}

%s3.5 #&#
\subsection{A superduality}\label{Ssupdu}

We have already collected all the necessary material to prove the lower bound
(\ref{lbd}) in Theorem~\ref{Tdecay}. Indeed, this follows from
Lemma~\ref{Lmonot}, which allows us to compare with a system
of coalescing random walks, for which the decay of the survival probability
and the density are given by Lemma~\ref{Lasmeet} and formula (\ref{ptau}).

The proof of the upper bound (\ref{ubd}) in Theorem~\ref{Tdecay} is more
involved. To prepare for this, in the present section, we will use
the graphical representation to construct (in terminology explained in
Section~\ref{Stermin}) a pathwise superdual to the cooperative
branching-coalescent.

Fix a graphical representation for the cooperative branching-coalescent, as
explained in Section~\ref{Sdefs}, consisting of Poisson point
processes $\lvec\omega(i),\rvec\omega(i)$ representing coalescing
jump events and
cooperative branching events, which occur on the whole time axis
${\mathbb R}$ (including
negative times). Ignoring the cooperative branching events for the moment
being, we define dual coalescing random walk arrows and \emph{dual
open paths} in
such a graphical representation as in Section~\ref{Scodu}.

Next, for any deterministic $s<u$, we define a \emph{dual $3$-path}
to be a triple of c\`adl\`ag functions $\gamma^k\dvtx [s,u)\to{\mathbb
Z}+\frac{1}{2}$ $(k=1,2,3)$
satisfying the following rules.\vspace*{6pt}

There exist times $s=t_0<\cdots<t_{n+1}=u$ such that:
\begin{longlist}[(2)]
\item[(1)] On each of the intervals $[t_{i-1},t_i)$ with $i=1,\ldots,n$
and on
$[t_n,u]$, the functions $\gamma^1,\gamma^2,\gamma^2$ are open dual
paths satisfying
$\gamma^1<\gamma^2<\gamma^3$.
\item[(2)] For each $t=t_i$ with $i\in1,\ldots,n$, one of the
following cases
occurs:
\begin{enumerate}[(a)]
\item[(a)]$t\in\lvec\omega(\gamma^k_t+\frac{1}{2})$ with $k=2,3$ (but
not $1$), and
$(\gamma^1_{t-},\gamma^2_{t-},\gamma^3_{t-})=(\gamma^k_t,\gamma
^k_t+1,\gamma^k_t+2)$,
\item[(b)]$t\in\rvec\omega(\gamma^k_t-\frac{1}{2})$ with $k=1,2$ (but
not $3$), and
$(\gamma^1_{t-},\gamma^2_{t-},\gamma^3_{t-})=(\gamma^k_t-2,\gamma
^k_t-1,\gamma^k_t)$.
\end{enumerate}
\end{longlist}
%

%
%f5 #&#
\begin{figure}%[h!]

\includegraphics{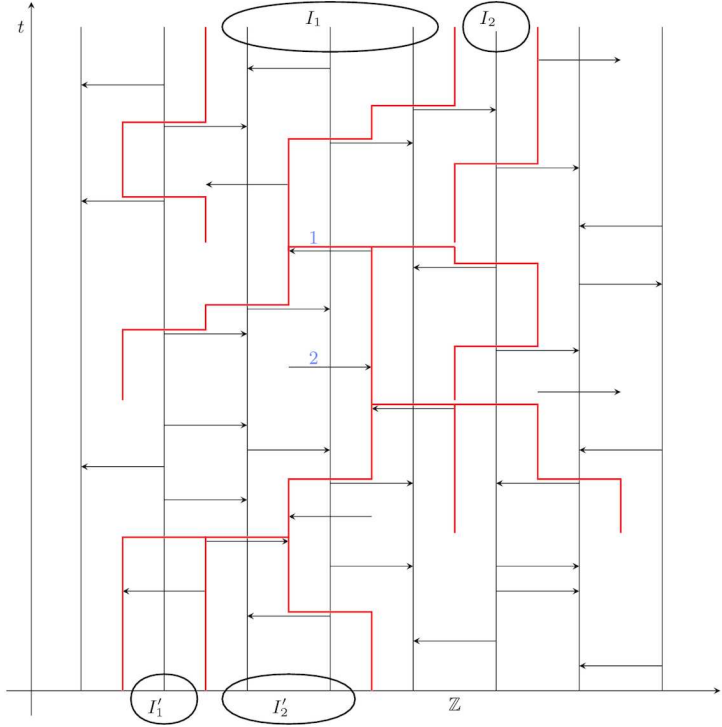}

%
%%need to make coordinates y1--y26
%%below 0.4 should be \y1
%
% -- (0.5,\vvw+0.05) -- (1.5,\vvw+0.05) -- (1.5,\vvu+0.1);
% -- (3.5,\vvz+0.05)-- (2.5,\vvz+0.05)-- (2.5,\vvu);
% -- (4.5,\vvu+0.1);
% -- (5.5,\vvp+0.05) -- (4.5,\vvp+0.05) -- (4.5,\vvm+0.1);
% -- (3.5,\vvi+0.05) -- (2.5,\vvi+0.05) -- (2.5,\vvc+0.05)
% -- (3.5,\vvc+0.05) -- (3.5,0);
% -- (2.5,\vvr+0.05)-- (1.5,\vvr+0.05) -- (1.5,\vvq+0.05)
% -- (0.5,\vvq+0.05) -- (0.5,\vvm+0.1);
% -- (4.5,\vvf+0.1);
% -- (5.5,\vvi+0.05) -- (6.5,\vvi+0.05) -- (6.5,\vvf+0.1);
% -- (1.5,0);
% -- (0.5,0);
%
%
%coordinates{(1.7,8) (3,8.3) (4.3,8) (3,7.7)};
%coordinates{(4.6,8) (5,8.3) (5.4,8) (5,7.7)};
%coordinates{(0.6,-0.1) (1,0.2) (1.4,-0.1) (1,-0.4)};
%coordinates{(1.7,-0.1) (2.5,0.2) (3.3,-0.1) (2.5,-0.4)};
%
%%\draw[->] (1.5,2)--(0.5,2);
%%\draw[->,very thick] (1,3.1)--(2,3.1);
%%\draw[->,very thick] (2.5,1)--(1.5,1);
%%\draw[->,very thick] (2,1.5)--(3,1.5);
%%\draw[->,very thick] (2.5,3.7)--(3.5,3.7);
%%\draw[->,very thick] (5,0.9)--(4,0.9);
%%\draw[->] (5,1.5)--(6,1.5);
%%\draw[->] (4.5,2.2)--(5.5,2.2);
%%\draw[->] (5,3)--(4,3);
%%
%%\draw[very thick] (1,1)--(1,3.1);
%%\draw[very thick] (2,0)--(2,1.5);
%%\draw[very thick] (2,3.1)--(2,4);
%%\draw[very thick] (3,0)--(3,4);
%%\draw[very thick] (4,0.9)--(4,4);
%%\draw[very thick] (5,0)--(5,0.9);
%%
\caption{Superduality: if $\eta_t\cap I_1\neq\varnothing$ and $\eta
_t\cap
I_2\neq\varnothing$, then there must exist a backward 3-path as drawn such
that $\eta_0\cap I'_1\neq\varnothing$ and $\eta_0\cap I'_2\neq
\varnothing$.
Between times when the 3-path renews itself, it consists of
three dual coalescing random walk paths which form the boundaries of two
adjacent intervals as in Figure~\protect\ref{figcodu}. Cooperative branching
arrows such as the one marked ``1'' \emph{may} be used to renew the 3-path
by splitting one of its paths into three new paths, but they do not
\emph{need} to be used such as the cooperative branching arrow 2.}
\label{figsubdu}
\end{figure}

An example of a dual 3-path is drawn in Figure~\ref{figsubdu}. In the
absense of cooperative branching events, the three paths $\gamma
^1,\gamma^2,\gamma^3$
evolve as dual coalescing random walk paths, which, however, are not
allowed to
coalesce (if the dual coalescing random walks coalesce then the 3-path ends).
If either $\gamma^2$ or $\gamma^3$ (but not $\gamma^1$) hits the
head of a
cooperative branching arrow pointing to the left, then we may forget
about the
three old paths and start anew with three new backward random walks
from the
positions $i,i+1$, and $i+2$, where $i \in{\mathbb Z}+ \frac{1}{2}$
is the location of the head of the
cooperative branching arrow. A similar rule applies for cooperative branching
arrows pointing to the right. We say that a dual 3-path \emph{renews}
itself at
such an instance. Note\vadjust{\goodbreak} that cooperative branching events \emph{may} be
used to
renew the dual 3-path, but they do not \emph{need} to be used. As a result,
there may be many different dual 3-paths starting from a given initial state
$(\gamma^1_u,\gamma^2_u,\gamma^3_u)$ and running backwards in time.
It is not hard to
see that the times when a dual 3-path renews itself can a.s. be read
off from
the path;\vspace*{2pt} that is, all information is contained in the triple of c\`
adl\`ag functions
$(\gamma^1,\gamma^2,\gamma^2)$.

Recall the notation introduced in (\ref{haint}). We let
%
%e82 #&#
\begin{equation}
\label{Xi2} \Xi_{+,2}:= \bigl\{\bigl(\langle i,j\rangle,\langle j,k
\rangle\bigr)\dvtx i,j,k\in{\mathbb Z}+\tfrac{1}{2}, i<j<k \bigr\}
\end{equation}
denote the space whose elements are pairs $(I_1,I_2)=(\langle
i,j\rangle,\langle j,k\rangle)$
of adjacent, discrete, nonempty, finite intervals in ${\mathbb Z}$.
The usefulness of
dual 3-paths lies in the following fact.

%le14 #&#
\begin{lemma}[(Dual 3-paths)]\label{L3path}
Let $(\eta_t)_{t\geq0}$ be a cooperative branching--coales\-cent
constructed
with a graphical representation as described in Section~\ref{Sdefs}. Let
$0\leq s<u$, and let $(I_1,I_2)\in\Xi_{+,2}$ be a pair of adjacent
intervals in
${\mathbb Z}$. Then, a.s. on the event
%
%e83 #&#
\begin{equation}
\label{inABu} \eta_u\cap I_1\neq\varnothing\quad
\mbox{and}\quad\eta_u\cap I_2\neq\varnothing,
\end{equation}
there exists a $(I'_1,I'_2)\in\Xi_{+,2}$ and a dual 3-path
$(\gamma^1_t,\gamma^2_t,\gamma^3_t)_{t\in[s,u]}$ with
%
%e84 #&#
\begin{equation}
(I_1,I_2)= \bigl(\bigl\langle\gamma^1_u,
\gamma^2_u\bigr\rangle,\bigl\langle\gamma
^2_u,\gamma^3_u\bigr\rangle
\bigr) \quad\mbox{and}\quad\bigl(I'_1,I'_2
\bigr)= \bigl(\bigl\langle\gamma^1_s,
\gamma^2_s\bigr\rangle,\bigl\langle\gamma^2_s,
\gamma^3_s\bigr\rangle\bigr),
\end{equation}
such that
%
%e85 #&#
\begin{equation}
\label{inABs} \eta_s\cap I'_1\neq
\varnothing\quad\mbox{and}\quad\eta_s\cap I'_2
\neq\varnothing.
\end{equation}
\end{lemma}

\begin{pf}
In the absence of cooperative branching
events, there exist unique
dual open paths $(\gamma^1_t,\gamma^2_t,\gamma^3_t)_{t\in[s,u]}$
starting at time $u$
from the boundaries of $I_1$~and~$I_2$, and these form a dual 3-path,
by the
definition of the latter, if and only if they do not coalesce until
time $s$.
Thus, in this case, coalescing random walk duality (\ref{walkdual})
tells us
that the events in (\ref{inABu}) and (\ref{inABs}) are in fact
a.s. equivalent.

In general, in the presence of cooperative branching events, let us define,
for $s\leq t\leq u$,
%
%e86 #&#
\begin{eqnarray}
\qquad {\mathcal J}_t&:=& \bigl\{
\bigl(I'_1,I'_2\bigr)\in\Xi
_{+,2}\dvtx \exists\mbox{ dual 3-path } \bigl(\gamma^1_s,
\gamma^2_s,\gamma^3_s
\bigr)_{s\in[t,u]}\mbox{ s.t.}
\nonumber\\[-8pt]\\[-8pt]
&&\hspace*{5pt} (I_1,I_2)= \bigl(\bigl\langle\gamma^1_u,
\gamma^2_u\bigr\rangle,\bigl\langle\gamma^2_u,
\gamma^3_u\bigr\rangle\bigr) \mbox{ and }
\bigl(I'_1,I'_2\bigr)= \bigl(
\bigl\langle\gamma^1_t,\gamma^2_t
\bigr\rangle,\bigl\langle\gamma^2_t,\gamma^3_t
\bigr\rangle\bigr) \bigr\}.\nonumber
\end{eqnarray}
It suffices to prove that if a cooperative branching event takes place
at time
$t$ and
%
%e87 #&#
\begin{equation}
\label{inII} \exists\bigl(I'_1,I'_2
\bigr)\in{\mathcal J}_t\qquad\mbox{s.t. } \eta_t\cap
I'_1\neq\varnothing\mbox{ and }
\eta_t\cap I'_2\neq\varnothing,
\end{equation}
then the same is true at time $t-$, that is, just before the cooperative
branching event. By symmetry, it suffices to consider the case
$t\in\lvec\omega(i)$ for some $i\in{\mathbb Z}$. By assumption,
(\ref{inII}) holds at time
$t$, so there exist $(I'_1,I'_2)\in{\mathcal J}_t$ with $\eta_t\cap
I'_1\neq\varnothing$
and $\eta_t\cap I'_2\neq\varnothing$. The only way in which this can
fail to
hold at time $t-$ is that the cooperative branching event has
introduced a
particle (at $i-1$)
into either $I'_1$ or $I'_2$, while this interval was
empty at time $t-$. For this to happen, the arrow associated with
$\lvec\omega(i)$
must point into $I'_1$ or $I'_2$ from the outside and $i$ and $i+1$
must both have been occupied by a particle at time $t-$. But
then, by the way dual 3-paths may renew themselves, we have
$(\{i\},\{i+1\})\in{\mathcal J}_{t-}$
and hence (\ref{inII}) is also
satisfied in this case.
\end{pf}

We claim that Lemma~\ref{L3path} actually gives rise to a Markov process
that, using terminology defined in Section~\ref{Stermin}, is a pathwise
superdual to the cooperative branching-coalescent. To see this,
we change the notation introduced in the proof of Lemma~\ref{L3path} slightly.
For any finite subset ${\mathcal J}_0\subset\Xi_{+,2}$ and fixed
$u\in{\mathbb R}$, define a
Markov process $({\mathcal J}_t)_{t\geq0}$ taking values in the finite
subsets of
$\Xi_{+,2}$, by
%
%e88 #&#
\begin{eqnarray}\label{supdu}
{\mathcal J}_t&:=& \bigl\{\bigl(I'_1,I'_2
\bigr)\in\Xi_{+,2}\dvtx  \exists(I_1,I_2)\in{\mathcal
J}_0\mbox{ and a dual 3-path}\nonumber
\\
&&\hspace*{5pt}\bigl(\gamma^1_s,\gamma^2_s,
\gamma^3_s\bigr)_{s\in
[u-t,u]}\mbox{ s.t. }
(I_1,I_2)= \bigl(\bigl\langle\gamma^1_u,
\gamma^2_u\bigr\rangle,\bigl\langle\gamma
^2_u,\gamma^3_u\bigr\rangle
\bigr)
\\
&&\hspace*{63pt} \mbox{and }\bigl(I'_1,I'_2
\bigr)= \bigl(\bigl\langle\gamma^1_{u-t},
\gamma^2_{u-t}\bigr\rangle,\bigl\langle\gamma
^2_{u-t},\gamma^3_{u-t}\bigr\rangle
\bigr) \bigr\}.\nonumber
\end{eqnarray}
Letting $\psi$ denote the duality function
%
%e89 #&#
\begin{equation}
\psi(\eta,{\mathcal J}):=1_{\{\exists(I_1,I_2)\in
{\mathcal J}\ \mathrm{s.t.}\ \eta\cap
I_1\neq\varnothing\ \mathrm{and}\ \eta\cap I_2\neq\varnothing\}},
\end{equation}
the proof of Lemma~\ref{L3path} then shows that the function
%
%e90 #&#
\begin{equation}
[0,u]\ni t\mapsto\psi(\eta_{t-},{\mathcal J}_{u-t})
\end{equation}
is a.s. nonincreasing; that is, the process $({\mathcal J}_t)_{t\geq
0}$ is a pathwise
superdual of $(\eta_t)_{t\geq0}$.

%s3.6 #&#
\subsection{Extinction of the superdual}\label{Ssupex}

In this section, we show that the superdual from (\ref{supdu}) dies out
a.s. (i.e., ${\mathcal J}_t=\varnothing$ eventually) if the cooperative
branching rate
satisfies $\lambda<1/2$. To keep the argument simple, and since this
is all we
will need in the end, we will only show this for the simplest possible initial
state, where ${\mathcal J}_0=\{(I_1,I_2)\}$ contains only a single pair
of adjacent
intervals, and these both have length one. We fix some $u\in{\mathbb
R}$. For each $t>0$
we consider the quantity
%
%e91 #&#
\begin{eqnarray}\label{Ndef}
N_t&:=&\mbox{the number of distinct dual 3-paths }
\bigl(\gamma^1_s,\gamma^2_s,
\gamma^3_s\bigr)_{s\in[u-t,u]}
\nonumber\\[-8pt]\\[-8pt]
&&\mbox{such that } \bigl(\gamma^1_u,
\gamma^2_u,\gamma^3_u\bigr)=
\bigl(-\tfrac{1}{2},\tfrac{1}{2},\tfrac{{3}}{{2}}\bigr).\nonumber
\end{eqnarray}
The next lemma not only shows that the probability that $N_t\neq0$ tends
to zero as $t\to\infty$, but more importantly also determines the right
speed of decay.

%le15 #&#
\begin{lemma}[(Expected number of dual 3-paths)]\label{L3mean}
Let $K$ be the constant from~(\ref{meet2}), and let $N_t$ be as in
(\ref{Ndef}). Then
%
%e92 #&#
\begin{equation}
\label{3mean} {\mathbb E}[N_t]\leq K \Biggl(\sum
_{n=1}^\infty(2\lambda)^nn^{5/2}
\Biggr)t^{-3/2}\qquad(t\geq0).
\end{equation}
\end{lemma}

\begin{pf}
Let $\tau_{\langle3\rangle}$ be the first
meeting time of three walkers
as in (\ref{tau23}), and set
%
%e93 #&#
\begin{equation}
G(t):={\mathbb P}[t\leq\tau_{\langle3\rangle}]\qquad(t\geq0).
\end{equation}
We may distinguish dual 3-paths according to how often they renew themselves
on the interval $[u-t,u]$. The probability that
there is a dual 3-path on $[u-t,u]$ that never renews itself is then
given by
$G(t)$ (recall that appropriate cooperative branching events may be
used for renewal but
that they do not have to be used). Since there are four ways in which a
path can renew itself, each of
which has rate $\lambda/2$, the probability that there is a dual
3-path on
$[u-t,u]$ that renews itself in the time interval $(u-s,u-s-\mathrm
{d}s)$ is
%
%e94 #&#
\begin{equation}
G(s)\cdot(2\lambda\,\mathrm{d}s)\cdot G(t-s).
\end{equation}
Thus the expected number of paths that renew themselves exactly once
during the
time interval $[u-t,t]$ is given by
%
%e95 #&#
\begin{equation}
2\lambda\int_0^t \mathrm{d}s\, G(s)G(t-s)=2
\lambda G\ast G(t),
\end{equation}
where $\ast$ denotes the convolution of two functions. Similarly, the expected
number of paths that renew themselves exactly $n$ times during the time
interval $[u-t,t]$ is given by
%
%e96 #&#
\begin{equation}
(2\lambda)^nG^{\ast n}(t),
\end{equation}
where $G^{\ast n}$ denotes the $n$th convolution power of $G$, and
hence
%
%e97 #&#
\begin{equation}
\label{EN} {\mathbb E}[N_t]=\sum_{n=1}^\infty(2
\lambda)^nG^{\ast n}(t).
\end{equation}

Let $G_1,G_2$ be functions satisfying
%
%e98 #&#
\begin{equation}
\label{Gicond} \int_0^\infty G_i(t)
\,\mathrm{d}t=1 \quad\mbox{and}\quad0 \leq G_i(t)\leq
K_it^{-\alpha}\qquad(i=1,2, t\geq0),
\end{equation}
and let $0<p<1$. Then
%
%e99 #&#
\begin{eqnarray}
 G_1\ast G_2(t)&=&\int_0^t
\mathrm{d}s\, G_1(s)G_2(t-s)
\nonumber
\\
&=&\int_{pt}^t\mathrm{d}s\,
G_1(s)G_2(t-s) +\int_{(1-p)t}^t
\mathrm{d}s\, G_2(s)G_1(t-s)
\nonumber
\\
&\leq&\int_{pt}^t\mathrm{d}s
K_1s^{-\alpha}G_2(t-s) +\int_{(1-p)t}^t
\mathrm{d}s K_2s^{-\alpha}G_1(t-s)
\\
&\leq& K_1(pt)^{-\alpha}\int_{pt}^t
\mathrm{d}s\, G_2(t-s) +K_2\bigl((1-p)t\bigr)^{-\alpha}
\int_{(1-p)t}^t\mathrm{d}s\, G_1(t-s)
\nonumber
\\
&\leq&\bigl(K_1p^{-\alpha}+K_2(1-p)^{-\alpha}
\bigr)t^{-\alpha},
\nonumber
\end{eqnarray}
where in the last step we have used that $G_1$ and $G_2$ have integral one.
By induction, we get for functions $G_1,\ldots,G_n$ the estimate
%
%e100 #&#
\begin{equation}
\label{convpow} G_1\ast\cdots\ast G_n(t) \leq
\bigl(K_1p_1^{-\alpha}+\cdots+K_np_n^{-\alpha}
\bigr)t^{-\alpha},
\end{equation}
where $p_1,\ldots,p_n$ are nonnegative numbers summing up to one.

In our case, condition (\ref{Gicond}) is satisfied by Corollary~\ref
{Lmeet2} and Lemma~\ref{Lmeanmeet} since
%
%e101 #&#
\begin{equation}
\int_0^\infty\mathrm{d}t G(t) =\int
_0^\infty\mathrm{d}t\, {\mathbb P}[t\leq
\tau_{\langle
3\rangle}]={\mathbb E}[\tau_{\langle3\rangle}]=1.
\end{equation}
Hence, setting $p_i=1/n$ and $\alpha= 3/2$ in (\ref{convpow}) we
obtain in our set-up the estimate
%
%e102 #&#
\begin{equation}
G^{\ast n}(t)\leq K\cdot n\cdot(1/n)^{-3/2}\cdot
t^{-3/2}=Kn^{5/2}t^{-3/2},
\end{equation}
which by (\ref{EN}) yields (\ref{3mean}).
\end{pf}

%s3.7 #&#
\subsection{Algebraic decay}

In this section we prove Theorem~\ref{Tdecay}. We start with some preparatory
lemmas. The first concerns an upper bound for the decay of the density of
pairs of particles: at least for $\lambda<1/2$ we obtain the same
rate of decay as in the case $\lambda=0$ without
cooperative branching; see (\ref{ptau}) and Corollary~\ref{Lmeet2}.

%le16 #&#
\begin{lemma}[(Density of pairs)]\label{Lpair}
Let $(\eta_t)_{t\geq0}$ be a cooperative branching-coalescent with cooperative
branching rate $\lambda<1/2$, started in an arbitrary initial law.
Let $K$ be the constant from (\ref{meet2}), and let
%
%e103 #&#
\begin{equation}
\label{Kac} K':=K\sum_{n=1}^\infty(2
\lambda)^nn^{5/2}<\infty.
\end{equation}
Then
%
%e104 #&#
\begin{equation}
{\mathbb P} \bigl[\{0,1\}\subset\eta_t \bigr]\leq
K't^{-3/2}\qquad(t\geq0).
\end{equation}
\end{lemma}

\begin{pf}
By\vspace*{1pt} Lemma~\ref{L3path}, the
probability that $\{0,1\}\subset\eta_t$ is bounded from above by the
probability that there exists a dual 3-path $(\gamma^1_s,\gamma
^2_s,\gamma^3_s)_{s\in[0,t]}$
with $(\gamma^1_t,\gamma^2_t,\gamma^3_t)=(-\frac{1}{2},\frac
{1}{2},\frac{{3}}{{2}})$ such that
%
%e105 #&#
\begin{equation}
\eta_0\cap\bigl\langle\gamma^1_0,
\gamma^2_0\bigr\rangle\neq\varnothing\quad\mbox{and}\quad
\eta_0\cap\bigl\langle\gamma^2_0,
\gamma^3_0\bigr\rangle\neq\varnothing.
\end{equation}
By Lemma~\ref{L3mean} we can estimate this from above, uniformly in the
initial law of $X$, by
%
%e106 #&#
\begin{equation}
{\mathbb P}[N_t>0]\leq{\mathbb E}[N_t]\leq
K't^{-3/2}\qquad(t\geq0).
\end{equation}\upqed
\end{pf}

%le17 #&#
\begin{lemma}[(Expected number of occupied pairs)]\label{LEpair}
Let $(\eta_t)_{t\geq0}$ be a cooperative branching-coalescent with
cooperative
branching rate $\lambda<1/2$, started in $\eta_0=\{0,1\}$, and let
$K'$ be the
constant from (\ref{Kac}).
Then
%
%e107 #&#
\begin{equation}
\label{Epair} {\mathbb E} \bigl[ \bigl| \bigl\{i\in{\mathbb Z}\dvtx \{i,i+1\}
\subset\eta
_t \bigr\} \bigr| \bigr] \leq K't^{-3/2}\qquad(t
\geq0).
\end{equation}
\end{lemma}

\begin{pf}
By Lemma~\ref{L3path}, for each $i\in
{\mathbb Z}$, the probability
%
%e108 #&#
\begin{equation}
{\mathbb P} \bigl[\{i,i+1\}\subset\eta_t \bigr]
\end{equation}
is\vspace*{1pt} bounded from above by the probability that there exists a dual 3-path
$(\gamma^1_s,\gamma^2_s,\gamma^3_s)_{s\in[0,t]}$ with
$(\gamma^1_t,\gamma^2_t,\gamma^3_t)=(i-\frac{1}{2},i+\frac
{1}{2},i+\frac{{3}}{{2}})$ such that
$\gamma^2_0=\frac{1}{2}$. By translation invariance, this is the same
as the
probability that there exists a dual 3-path with
$(\gamma^1_t,\gamma^2_t,\gamma^3_t)=(-\frac{1}{2},\frac
{1}{2},\frac{{3}}{{2}})$ such that
$\gamma^2_0=-i+\frac{1}{2}$. Summing over all $i\in{\mathbb Z}$ and using
Lemma~\ref{L3mean}, this implies that
%
%e109 #&#
\begin{equation}
{\mathbb E} \bigl[ \bigl| \bigl\{i\in{\mathbb Z}\dvtx \{i,i+1\}\subset\eta_t
\bigr\}  \bigr| \bigr] \leq{\mathbb E}[N_t]\leq K't^{-3/2}
\qquad(t\geq0).
\end{equation}\upqed
\end{pf}

\begin{pf*}{Proof of Theorem~\ref{Tdecay}} By Lemma~\ref
{Lmonot}, we may
stochastically bound an arbitrary cooperative branching-coalescent by a
cooperative branching-coalescent with $\lambda=0$, that is, a
system of coalescing
random walks. Thus it suffices to prove the lower bound in (\ref{lbd}) only for
$\lambda=0$. For such systems, using notation introduced in~(\ref{tau23}), we have
that
%
%e110 #&#
\begin{equation}
{\mathbb P} \bigl[\bigl|\eta^{\{0,1\}}_t\bigr|\geq2 \bigr]={\mathbb P}[t
\leq\tau_{\langle2\rangle}] \quad\mbox{and}\quad{\mathbb P}\bigl[0\in
\eta^{\mathbb Z}_t\bigr]={\mathbb P}[t\leq\tau_{\langle
2\rangle}],
\end{equation}
where the second equality is (\ref{ptau}). By Lemma~\ref{Lasmeet}, there
exists a constant $c>0$ such that
%
%e111 #&#
\begin{equation}
{\mathbb P}[t\leq\tau_{\langle2\rangle}]\geq c t^{-1/2}\qquad(t\geq0).
\end{equation}
This completes the proof of the lower bound in (\ref{lbd}).

To get also the upper bound in (\ref{ubd}), define $p_t(\cdots)$ as in
(\ref{ptdef}) for the process $\eta^{\mathbb Z}$. Since, by
Theorem~\ref{Tphase}(a),
$p_t(1)\to0$ as $t\to\infty$, formula (\ref{difp}) tells us that
%
%e112 #&#
\begin{eqnarray}
p_t(1) &=& -\int_t^\infty\mathrm{d}s
\frac{\partial}{\partial{s}} p_s(1)\nonumber
\\
&=& \int_t^\infty
\mathrm{d}s \bigl((1-\lambda)p_s(11)+\lambda p_s(111)
\bigr)
\\
&\leq& \int_t^\infty\mathrm{d}s\, p_s(11).\nonumber
\end{eqnarray}
Since $p_s(11)\leq K's^{-3/2}$ by Lemma~\ref{Lpair}, we find that
%
%e113 #&#
\begin{equation}
{\mathbb P}\bigl[0\in\eta^{\mathbb Z}_t\bigr]\leq
K'\int_t^\infty\mathrm{d}s\,
s^{-3/2} =2 K't^{-1/2}.
\end{equation}
Similarly, the indicator function on the event $\{|\eta^{\{0,1\}
}_t|\geq2\}$ decreases at rate $1$ whenever
$\eta^{\{0,1\}}_t=\{i,i+1\}$ for some $i \in{\mathbb Z}$ due to an
appropriate random walk
step and subsequent coalescence. Thus, by Lemma~\ref{LEpair},
%
%e114 #&#
\begin{eqnarray}
-\frac{\partial}{\partial{t}}{\mathbb P} \bigl[\bigl|\eta^{\{0,1\}
}_t\bigr|\geq2
\bigr] &=&{\mathbb P} \bigl[\eta^{\{0,1\}}_t =\{i,i+1\}\mbox{ for
some } i \in{\mathbb Z} \bigr]\nonumber
\\
&\leq&
{\mathbb P} \bigl[ \{i,i+1\} \subset\eta^{\{0,1\}}_t
\mbox{ for some } i \in{\mathbb Z} \bigr]
\\
&\leq&
K't^{-3/2}\qquad(t\geq0).\nonumber
\end{eqnarray}
Hence, using Theorem~\ref{Tphase}(b) we find that
%
%e115 #&#
\begin{equation}
{\mathbb P} \bigl[\bigl|\eta^{\{0,1\}}_t\bigr|\geq2 \bigr]\leq\int
_t^\infty\mathrm{d}s\, K's^{-3/2}
=2K't^{-1/2}\qquad(t\geq0).
\end{equation}\upqed
\end{pf*}

\section*{Acknowledgements}
We thank Wolfgang K\"onig, Jochen Blath and\break Rongfeng Sun for useful
discussions. %This research was in part supported by the DFG through priority programme 1590.

%suskaldyti doi

% imsref loaded by linak, 2014-05-26 15:56:29

\printaddresses
\end{document}